\documentclass[a4paper]{article}
\usepackage{amsmath}
    \usepackage{amsxtra}
    \usepackage{amstext}
    \usepackage{amssymb}
    \usepackage{latexsym}
\usepackage[latin1]{inputenc}
\usepackage{setspace}
\input xy

\xyoption{all}

\xyoption{poly}

\newenvironment{changemargin}[2]{%
  \begin{list}{}{%
    \setlength{\topsep}{0pt}%
    \setlength{\leftmargin}{#1}%
    \setlength{\rightmargin}{#2}%
    \setlength{\listparindent}{\parindent}%
    \setlength{\itemindent}{\parindent}%
    \setlength{\parsep}{\parskip}%
  }%
  \item[]}{\end{list}}

\newcommand{\nin}{\mbox{$\in \!\!\!\!\!/\,$}}

\newtheorem{prop}{Proposition}[section]
\newtheorem{thm}{Theorem}[section]
\newtheorem{thm'}{Theorem}
\newtheorem{cor}{Corollary}[section]
\newtheorem{lem}{Lemma}[section]

 \newcommand{\calF}{\mathcal F}

\newcommand{\mA}{\mathbb A}
 \newcommand{\mC}{\mathbb C} 
 \newcommand{\mF}{\mathbb F}

  \newcommand{\mP}{\mathbb P}
\newcommand{\mQ}{\mathbb Q} \newcommand{\mR}{\mathbb R}

\newcommand{\mZ}{\mathbb Z}

\newcommand{\gra}{\alpha} \newcommand{\grb}{\beta}       
  
      \newcommand{\grs}{\sigma}
\newcommand{\gro}{\omega}      
\newcommand{\grt} {\theta}

 \newcommand{\grD}{\Delta}

\newcommand{\vuoto}      {\varnothing}

\newcommand{\cech}       {\vee}

\newcommand{\ol}         {\overline}
\newcommand{\wt}         {\widetilde}
\newcommand{\ogni}         {\forall}

\begin{document}

\begin{center}
\bf \large FANO SYMMETRIC VARIETIES WITH LOW RANK
\end{center}

\

\begin{center}
  \large   Alessandro Ruzzi
\end{center}

\

\begin{abstract}
The symmetric projective varieties of rank one are all smooth and Fano by a classic result of Akhiezer.
We classify the locally factorial (respectively smooth) projective symmetric $G$-varieties of rank 2 which are Fano.  When $G$ is semisimple we classify also the locally factorial (respectively smooth) projective symmetric $G$-varieties of rank 2 which are only quasi-Fano. Moreover, we classify the Fano symmetric $G$-varieties of rank 3   obtainable from a wonderful variety by a sequence of blow-ups along $G$-stable varieties. Finally, we classify the Fano symmetric  varieties of arbitrary rank which are obtainable from a wonderful variety  by a sequence of blow-ups along closed orbits.
\end{abstract}

{\small\ \ \  \textit{keywords}: Symmetric varieties, Fano
varieties.}


\

A Gorenstein (projective) normal algebraic variety $X$ over $\mathbb{C}$ is
called a Fano variety if   the anticanonical divisor is ample.  The
Fano surfaces are classically called  Del Pezzo surfaces. The
importance of Fano varieties in the theory of higher dimensional
varieties is similar to the importance of Del Pezzo surfaces in the
theory of surfaces. Moreover  Mori's program predicts that every
uniruled variety is birational to a fiberspace whose general fiber
is a Fano variety (with terminal singularities).

Let $\grt$ be an involution of a reductive group $G$ (over $\mC$) and let $H$ be a closed subgroup of $G$ such that  $G^\grt\subset H\subset N_G(G^\grt)$.  A symmetric variety is a normal $G$-variety with an open orbit isomorphic to $G/H$. The symmetric varieties are a generalization of the toric varieties. The toric smooth Fano varieties with rank at most four are been classified. By \cite{AlBr}, Theorem 4.2  there is only a finite number of Fano smooth
symmetric varieties with a fixed open orbit. In
\cite{Ru1} we have classified the smooth compact symmetric varieties
with Picard number one and $G$ semisimple, while in \cite{Ru2} we have given an explicitly
geometrical description of such varieties; they  are automatically Fano.

In this work, we want to classify the Fano  symmetric varieties  with
low rank (and $G$ semisimple).  First, we consider a special case of arbitrary rank.
We say that a variety $X$ is quasi $\mQ$-Fano if   $-K_{X}$ is a nef and big $\mathbb{Q}$-divisor.
Fixed an open orbit $G/H$ with $G$ semisimple, there is a unique maximal compactification between the ones which have  only one closed orbit. Such variety   is called the standard compactification. If it is also smooth, it is called the wonderful compactification; this is the case, for example, if $H=N_G(G^\grt)$ (see \cite{dCP1}, Theorem 3.1).
We prove that the standard symmetric varieties are all quasi $\mQ$-Fano and we describe when they are Fano. We determine also the symmetric Fano varieties obtainable from a wonderful one by a sequence of blow-ups along closed orbits. In
particular, we prove that   such a variety must be either a wonderful
one or the blow-up of a wonderful one along the unique closed orbit.

Next we consider the symmetric varieties of rank at most three. The rank of a symmetric variety $X$ is defined as  the rank of $\mC(X)^{(B)}/\mC^*$, where $B$ is any fixed Borel subgroup of $G$. The symmetric varieties with rank one are all wonderful; moreover one can
show that they are isomorphic, under the action of $Aut^{0}(X)$,
either to  a projective homogeneous variety $G/P$ with $P$
maximal, or to $\mathbb{P}^{n}\times \mathbb{P}^{n}$ (see \cite{A83}). Thus they are
all Fano.

We classify all the locally factorial (resp. smooth)  Fano symmetric varieties of
rank 2. When $G$ is semisimple, we classify also the locally factorial (resp. smooth)  symmetric varieties which are  only quasi-Fano.  In the proof of such result we obtain a classification of the toroidal Fano varieties of rank 2 with $G$ semisimple (without assumption on the regularity).

Finally, we classify the smooth Fano symmetric varieties of rank three which
are obtainable from a wonderful one by a sequence of blow-ups along
$G$-subvarieties (in particular $G$ is semisimple). This class of varieties is quite large; indeed
any compact symmetric variety is dominated by a variety obtained from the
wonderful one by a sequence of blow-ups along $G$-subvarieties of
codimension two (see \cite{dCP2}, Theorem 2.4). This   result on 3-rank varieties can be generalized to varieties obtainable from a generic wonderful varieties of rank 3 by a sequence of blow-ups along $G$-subvarieties (without suppose $G/H$ symmetric).

\section{Introduction and notations}\label{intro}

In this section we  introduce the necessary  notations. The reader
interested to the embedding theory of spherical varieties can see \cite{Kn},
\cite{Br97} or \cite{T06}. In \cite{V90} is explained such theory in
the particular case of  symmetric varieties.

\subsection{First definitions}\label{1def}

Let $G$ be a connected reductive algebraic group over $\mathbb{C}$
and let $\theta$ be an involution of $G$. Given  a closed
subgroup $H$  such that $G^{\theta}\subset H\subset
N_{G}(G^{\theta})$, we say that   $G/H$ is a {\em symmetric space}
and that $H$ is a {\em symmetric subgroup}.
A normal $G$-variety is called a {\em spherical
variety} if it contains a dense $B$-orbit ($B$ is a  chosen Borel subgroup of $G$). We
say that a subtorus of $G$ is split if $\theta(t)=t^{-1}$ for all
its elements $t$; moreover it is a {\em maximal split torus} if it has
maximal dimension. A   maximal torus containing a
maximal split torus is {\em maximally split};   any
maximally split torus is $\theta$ stable (see \cite{T06}, Lemma
26.5). We fix arbitrarily a maximal split torus $T^{1}$ and a
maximally split torus  $T$ containing $T^{1}$. Let $R_{G}$ be the
root system of $G$ w.r.t. $T$. We can choose a Borel
subgroup $T\subset B$  such that, for any positive root
$\alpha$, either $\theta(\alpha)=\alpha$ or $\theta(\alpha)$ is
negative. Moreover, $BH$ is dense in $G$ (see \cite{dCP1}, Lemma 1.2
and Proposition 1.3). In particular,  every
normal equivariant embedding of $G/H$ is  spherical; we call it a {\em symmetric variety}.

We can assume that $G$ is the direct product of a simply
connected, semisimple group with a central split torus.

\subsection{Colored fans}\label{colored fans}

Now,  we  introduce some details about the classification
of the symmetric varieties by their colored fans (this
classification is defined more generally for spherical varieties).
Let $D(G/H)$ be the set of $B$-stable prime divisors of $G/H$; its
elements are called {\em colors}.   We say that a spherical variety is {\em simple} if it contains
one closed orbit. Let $X$ be a simple symmetric variety with closed
orbit $Y$. We define the set of \textit{colors of $X$} as the subset $D(X)$
of $D(G/H)$ consisting of the colors whose closure in $X$ contains
$Y$.  To each prime divisor $D$ of $X$, we can associate the
normalized discrete valuation $v_{D}$ of $\mathbb{C}(G/H)$ whose
ring is $\mathcal{O}_{X,D}$; $D$
is $G$-stable if and only if $v_{D}$ is $G$-invariant. Let $N$ be the  set of all $G$-invariant
valuations of $\mathbb{C}(G/H)$ taking value in $\mathbb{Z}$ and let
$N(X)$ be the set of the valuations associated to the $G$-stable
prime divisors of $X$. Observe that each irreducible component  of
$X\, \backslash\, (G/H)$ has codimension one, because $G/H$ is
affine. Let $S:=T/\,T\cap H\simeq T\cdot (eH/H)$. One can show that
the group $\mathbb{C}(G/H)^{(B)}/\,\mathbb{C}^{*}$ is isomorphic to
the character group $\chi(S)$ of $S$ (see \cite{V90}, \S2.3); in
particular, it is a free abelian group. We define the {\em rank} $l$ of $G/H$
as the rank of $\chi(S)$. We can identify the dual group
$Hom_{\mathbb{Z}}(\mathbb{C}(G/H)^{(B)}/\,\mathbb{C}^{*},\mathbb{Z})$
with the group  $\chi_{*}(S)$ of one-parameter subgroups of $S$; so
we can identify $\chi_{*}(S)_{ \mathbb{R}}$ with
$Hom_{\mathbb{Z}}(\chi(S),\mathbb{R})$. The  restriction map to
$\mathbb{C}(G/H)^{(B)}/\mathbb{C}^{*}$ is injective over $N$ (see
\cite{Br97}, \S3.1 Corollaire 3), so we can identify $N$ with a
subset of $\chi_{*}(S)_{ \mathbb{R}}$. We say  that $N$ is the
{\em valuation monoid} of $G/H$. For each color $D$, we define
$\rho(D)$ as the restriction of $v_{D}$ to $\chi(S)$. In general,
the map $\rho:D(G/H)\rightarrow\chi_{*}(S)_{ \mathbb{R}}$ is not
injective. Let $C(X)$ be the cone in $\chi_{*}(S)_{ \mathbb{R}}$ generated by
$N(X)$ and $\rho(D(X))$. We say that the pair $(C(X),D(X))$ is the {\em colored cone} of X; it determines univocally $X$ (see \cite{Br97}, \S3.3 Théorème).

Let $Y$ be an orbit of a symmetric variety  $X$. The set $\{x\in X\
|\  \overline{G\cdot x}\supset Y\}$ is an open simple $G$-subvariety
of $X$ with closed orbit $Y$, because any spherical variety contains only
a finite number of $G$-orbits. Let $\{X_{i}\}$ be the set of open
simple subvarieties of $X$ and  define the set of \textit{colors of $X$},
$D(X)$, as $\bigcup_{i\in I}D(X_{i})$. The family
$\calF(X):=\{(C(X_{i}),D(X_{i}))\}_{i\in I}$ is \vspace{0.2 mm} called the
{\em colored fan} of $X$ and determines completely $X$ (see \cite{Br97},
\S 3.4 Th\'{e}or\`{e}me 1). Moreover $X$ is compact if and only if
$cone(N)$ is contained in the support $|\calF(X)|:=\bigcup_{i\in I}C(X_{i}) $ of
the colored fan (see \cite{Br97}, \S 3.4 Th\'{e}orem\`{e} 2).

Given a symmetric variety $X$ we denote by $\Delta$ (or by $\grD_X$) the fan associated to the colored fan of $X$, by $\Delta(i)$ the set of $i$-dimensional cones in $\Delta$ and by $\Delta[p]$ the set of primitive generators of the 1-dimensional cones of $\Delta$. The fan  $\Delta$ is formed by all the faces of the cones $C$ such that there is a colored cone $(C,F)\in \mathcal{F}(X)$.
The toric varieties are a
special case of symmetric varieties.  If $X$ is a toric variety,
then $D(G/H)$ is empty and we need only to consider the fan  $\grD_X$
(actually the theory of colored fans
is a generalization of the classification of toric varieties by
fans).

\subsection{Restricted root system}\label{restricted root system}

To describe  the sets $N$ and $\rho(D(G/H))$, we associate a
root system  to $G/H$. We can identify $\chi(T^{1} )_{\mathbb{R}}$
with $\chi(S)_{\mathbb{R}}$ because  $[\chi(S):\chi(T^{1})]$ is finite. We call again $\theta$ the involution induced on
$\chi(T)_{\mathbb{R}}$. The inclusion $T^{1}\subset T$ induces an
isomorphism of $\chi(T^{1})_{\mathbb{R}}$ with the $(-1)$-eigenspace
of $\chi(T)_{\mathbb{R}}$ under the action of $\theta$ (see
\cite{T06}, \S 26). Denote  by $W_{G}$ the Weyl group of $G$ (w.r.t. $T$). We can identify $\chi(T^{1})_{\mathbb{R}}$ with its
dual $\chi_{*}(T^{1})_{\mathbb{R}}$ by the restriction $(\, \cdot
,\cdot)$ to $\chi(T^{1})_{\mathbb{R}}$  of a fixed $W_G$-invariant non-degenerate symmetric bilinear form on
$\chi(T)_{\mathbb{R}}$ . Let
$R^{0}_{G}$ be the set of roots fixed by $\theta$ and let
$R^{1}_{G}$ be  $R_{G} \setminus  R^{0}_{G}$. Let $R_G^{i,+}:=R_G^{i}\cap R_G^{+}$.

The set $R_{G,\theta}:= \{\beta-\theta(\beta)\ |\ \beta\in
R_{G}^{1}\}$ is a root system in
$\chi(S)_{\mathbb{R}}$ (see \cite{V90}, \S 2.3 Lemme), which we call
the {\em restricted root system} of $(G,\theta)$; we call the non zero
$\beta-\theta(\beta)$ the {\em restricted roots}. Usually we denote by
$\beta$ (resp.  by $\alpha$)  a root of $R_{G}$ (resp.
of $R_{G,\theta}$); often we denote by $\varpi$ (resp. by
$\omega$)  a weight of $R_{G}$ (resp. of $R_{G,\theta}$). We
denote by $\overline{R}_{G}=\{\beta_{1},...,\beta_{n}\}$ the basis
of $R_{G}$ associated to $B$ and by $\varpi_{1},...,\varpi_{n}$ the
fundamental weights of $R_{G}$. Let $\overline{R}_{G}^{i}$ be
$\overline{R}_{G}\cap R_{G}^{i}$. There is a permutation $\overline{\grt}$ of $\overline{R}^1_G$ such that, $\forall\grb\in \overline{R}^1_G$, $\grt(\grb)+\overline{\grt}(\grb)$ is a linear combination of roots in $\overline{R}^0_G$.  We denote by $\alpha_{1},...,\alpha_{s}$ the elements
of the basis $\overline{R}_{G,\theta}:=\{\beta-\theta(\beta)\, |\,
\beta\in \ol{R}_{G}^1\}$   of $R_{G,\theta}$. If $R_{G,\theta}$ is irreducible we order $\overline{R}_{G,\theta}$ as in \cite{Bou}.
Let $b_{i}$ be equal to $\frac{1}{2}$ if $2\alpha_{i}$ belongs to
$R_{G,\theta}$ and equal to one otherwise; for each $i$ we define
$\alpha_{i}^{\vee}$ as the coroot
$\frac{2b_{i}}{(\alpha_{i},\alpha_{i})}\alpha_{i}$. The set
$\{\alpha^{\vee}_{1} ,...,\alpha^{\vee}_{s}\}$ is a basis of the
dual root system $R^{\vee}_{G,\theta}$. We call the
elements of $R^{\vee}_{G,\theta}$ the restricted coroots.  Let
$\omega_{1},...,\omega_{s}$  be the fundamental weights of $R
_{G,\theta}$  w.r.t. $\{\alpha_{1} ,...,\alpha_{s}\}$ and
let $\omega^{\vee}_{1},...,\omega^{\vee}_{s}$  be the fundamental
weights of $R^{\vee}_{G,\theta}$ w.r.t.
$\{\alpha^{\vee}_{1} ,...,\alpha^{\vee}_{s}\}$. Let $C^{+}$ be the
positive closed Weyl chamber of $\chi(S)_\mathbb{R}$ and let $C^-:=-C^{+}$.  

We say that a dominant weight $\varpi\in \chi(T)$ is a {\em spherical
weight} if $V(\varpi)$ contains a non-zero vector fixed by
$G^{\theta}$. In this case, $V(\varpi)^{G^{\theta}}$ is one-dimensional and
$\theta(\varpi)=-\varpi$, so  $\varpi$ belongs to
$\chi(S)_{\mathbb{R}}$. One can show that set of dominant weights of
$R_{G,\theta}$ is the set of spherical weights and that $C^{+}$ is
the intersection of $\chi(S)_{\mathbb{R}}$ with the positive closed
Weyl chamber of $R_{G}$. Suppose
$\beta_{j}-\theta(\beta_{j})=\alpha_{i}$, then $\omega_{i}$ is a
positive multiple of $\varpi_{j}+\varpi_{\overline{\theta}(j)}$.
More precisely, we have $\omega_{i}=\varpi_{j}+\varpi_{\theta(j)}$
if $\overline{\theta}(j)\neq j$, $\omega_{i}=2\varpi_{j}$ if
$\overline{\theta}(j)= j$ and $\beta_{j}$ is orthogonal to
$R_{G}^{0}$ and $\omega_{i}=\varpi_{j}$ otherwise (see \cite{CM},
Theorem 2.3 or \cite{T06}, Proposition 26.4). We say that a spherical weight is \textit{regular} if it is strictly dominant as weight of the restricted root system.

\subsection{The sets $N$ and $D(G/H)$}\label{n+dg/h}

The set $N$ is equal to $C^{-}\cap\,\chi_{*}(S)$; in particular, it
consists of the lattice vectors of the rational, polyhedral, convex
cone $C^{-}=cone(N)$. The set $\rho(D(G/H))$ is equal to
$\ol{R}^\cech_{G,\grt}$ and any fibre $\rho^{-1}(\gra^\cech)$ contains  at most 2 colors.  For any simple spherical variety $X$, $N(X)$ is formed by the primitive generators of the  1-faces of $C(X)$ which are contained in $cone(N)$. When $X$ is symmetric, also $\rho(X)$ can be recovered by $C(X)$: its elements generate the 1-faces of $C(X)$ which are not contained in $C^-$. We say that $(G,\theta)$
indecomposable if the unique  normal, connected, $\theta$-stable
subgroup of $G$ is the trivial one. In this case the number of
colors  is at most equal to $rank\,(G/H)+1$. If $\grt$ is indecomposable then there are three possibilities: i) $G$ is simple; ii) $G=\dot{G}\times \dot{G}$ with $\dot{G}$ simple and $\theta(x,y)=(y,x)$; iii) $G=\mC^*$ and $\grt(t)=t^{-1}$. See \cite{W}, \S1.1 for a classification of the involution of a simple group. In \cite{W} (and in \cite{Ru2}, \S1) are also indicated the Satake diagrams of the indecomposable involutions. The Satake diagram of any involution $(G,\grt)$ is obtained from the Dynkin diagram of $G$ as follows: 1) the vertices corresponding to element of $\overline{R}^0_G$  (resp. of $\overline{R}^1_G$) are black (resp. white); 2) two simple roots $\grb_1,\grb_2\in  \overline{R}^1_G$ such that $\overline{\grt}(\grb_1)=\grb_2$ are linked by a double-headed arrow.

If $\sharp D(G/H)>rank(G/H)$ and $(G,\theta)$ is indecomposable, we have
two possibilities: 1) $G^{\theta}=H=N_{G}(G^{\theta})$; 2) $H=G^{\theta}$ and $[G^\grt: N_{G}(G^{\theta})]=2$. In
the last case any element of $N_{G}(G^{\theta})\backslash\,
G^{\theta}$ exchanges  two colors and $R_{G,\grt}$ has type $A_1$, $B_2$ or $C_n$. We say that a simple restricted root $\gra$ is exceptional if $\sharp\rho^{-1}(\gra^\cech)=2$  and $2\gra$ is a restricted root. In this case the irreducible factor of $R_{G,\grt}$
containing $\gra$ is associated to an indecomposable factor of $G/G^\grt$ as in 1). We say that also  $(G,\grt)$ and any symmetric variety (with open orbit $G/H$) are exceptional. We denote by $D_{\alpha}$
the sum of the colors in $\rho^{-1}(\alpha^{\vee})$ and by $D_\gro$ the $G$-stable divisor corresponding to $(\mR^{\geq0}\gro,\vuoto)\in\mF(X)$.

\subsection{Toroidal symmetric varieties}\label{toroidal}

In this section we want to define a special
class of varieties. We say that a spherical variety is {\em toroidal}
if $D(X)=\vuoto$. 
\textit{Let $(C,F)$ be a colored cone of $X$, we say that the blow-up of $X$ along the subvariety associated to $(C,F)$ is the blow-up of $X$ along $(C,F)$.} In the following of this section we suppose $G$ semisimple.
Then there is a special simple compactification of $G/H$ because $N_{G}(H)/H$ is finite. This compactification, called
the standard compactification $X_{0}$, is
associated to $(cone(N),\vuoto)$ and it is the maximal simple
compactification  of $G/H$ in the dominant order. We define $e_i$ as the primitive positive multiple of
$-\gro_i^\cech$ (in $\chi_*(S)$), so $\grD_{X_0}[p]=\{e_1,...,e_l\}$. The standard compactification is wonderful (i.e. it is also smooth) if and only if $\chi_*(S)=\bigoplus\mathbb{Z}e_i$.
De Concini and Procesi have proved that $X_0$ is wonderful if $H=N_G(G^\grt)$, or equivalently $\chi_*(S)=\bigoplus\mathbb{Z}\gro^\cech_i$ (see \cite{dCP1} Theorem 3.1).

The standard compactification $X_{0}$ contains an affine toric $S$-variety $Z_{0}$, which is a
quotient of an affine space by a finite group. The toroidal
varieties are the symmetric varieties which dominates the standard
compactification and are in one-to-one correspondence with the $S$-toric
varieties which dominates $Z_{0}$.

Let $P$ be the stabilizer of the $B$-stable affine  open set $U:=X_{0}\setminus
\bigcup_{D(G/H)}\overline{D}$. This open set is $P$-isomorphic to $R_{u}P\times Z_{0}$, where
$R_{u}P=\prod_{\beta \in R^{1,+}_G}U_{\beta}$ is
the unipotent radical of $P$ and $dim\, Z_{0}=rank\, X_{0}$. To any
toroidal variety $X$ we associated the inverse image $Z$ of $Z_{0}$
by the projection $X\rightarrow X_{0}$. Moreover,  $X\setminus\
\bigcup_{D(G/H)}\overline{D}$ is $P$-isomorphic to $R_{u}P\times Z$.
The toroidal varieties are also in one-to-one correspondence with a
class of compact toric varieties in the following way. To a symmetric
variety   $X$, we associate the closure $Z^{c}$ of $Z$ in $X$;
$Z^{c}$ is also the inverse image of $Z_{0}^{c}$. The fan of $Z$ is
the fan $\grD_X$ associated to the colored fan of $X$, while the fan of
$Z^{c}$ consists of the translates of the cones of $Z$ by the Weyl
group $W_{G,\theta}\cong N_{G^\grt}(T^1)/C_{G^\grt}(T^1)$ of $R_{G,\theta}$.


\subsection{The Picard group}\label{sez Picard}

The class group of a symmetric variety is generated by the classes
of the $B$-stable prime divisors modulo the relations $div(f)$ with
$f\in \mathbb{C}(G/H)^{(B)}$. Indeed $Cl(BH/H)=Pic(BH/H)$ is
trivial. Given $\omega\in \chi(S)$ we denote by $f_{\omega}$  the
element of $\mathbb{C}(G/H)^{(B)}$ with weight $\omega$ and such
that $f_{\omega}(H/H)=1$.

A Weyl divisor $\sum_{D\in D(G/H)} a_{D}D+       \sum_{E\in
N(X)}b_{E}E$ is a Cartier divisor if and only if, for any $(C,F)\in\calF(X)$, there is $h_{C}\in \chi(S)$ such that $h_{C}(E)=a_{E}$
$\ogni E\in C$ and $h_{C}(\rho(D))=a_{D}$    $\ogni D\in F$.
Let $PL'(X)$ be the set of functions on the support $|\calF(X)| $  such that: 1) are linear on each colored
cone; 2) are integer on $\chi_{*}(S)\cap |\calF(X)| $. Let
$L(X)\subset PL'(X)$ be the subset composed by the restrictions of
linear functions and let $PL(X):=PL'(X)/L(X)$. The
$\{h_{C}\}$, corresponding to any  Cartier divisor,   defines an element of $PL(X)$.  If $X$ is compact,
there is an exact sequence (see \cite{Br89}, Théorème 3.1):

\[0\rightarrow \bigoplus_{D\in D(G/H)\setminus D(X)}\mathbb{Z}D\rightarrow Pic(X)\rightarrow PL(X)\rightarrow0.\]

A Cartier divisor is globally generated (resp. ample) if and only if
the associated function is convex (resp. strictly convex) and
$h_{C}(\rho(D))\leq a_{D}$ (resp. $h_{C}(\rho(D))< a_{D}$) $\ogni (C,F)\in\calF(X)$ and $\ogni D\in D(G/H)\setminus\, F$. Given any
linearized line bundle $L$, the space $H^{0}(X,L)$ is a multiplicity free $G$-module and,
if $L$ is globally generated, the highest weights of $H^{0}(X,L)$
are the elements of $\chi(S)\cap hull(\{h_{C}\}_{dim\,C=l})$, where $hull(\{x_1,...,x_m\})$ is the convex hull of  $x_1,...,x_m$ (see \cite{Br89},
\S3). Thus, a Cartier divisor on a projective symmetric variety is nef if and only if it is globally generated. Moreover, a nef $G$-stable  Cartier divisor on a projective symmetric variety is big if and only if the associated piecewise linear function $h$ is such that $(\sum_{C\in\grD(l)}h_C,R^\cech)\neq0$ for each irreducible factor $R^\cech$ of $R^\cech_{G,\grt}$ (see \cite{Ru09} Theorem 4.2). In particular, when $\grt$ is indecomposable every non-zero nef $G$-stable divisor is big. When $X$ is toroidal we have an exact split sequence

\[0\rightarrow Pic(X_{0})\rightarrow Pix(X)\rightarrow Pic(Z)\rightarrow0.\]

A normal    variety $X$ is locally factorial if the Picard group is isomorphic to the class group, while $X$ is $\mathbb{Q}$-factorial if $Pic(X)_{\mQ}\cong Cl(X)_{\mQ}$. A simple symmetric variety associated to a colored cone $(C,F)$ is locally factorial if: i) $C$ is generated by a subset of a basis of $\chi_*(S)$ and ii) $\rho$ is injective over $F$ (see \cite{Br97b} for a general statement in the spherical case). When the variety is toroidal the locally factoriality is equivalent to the smoothness.

An anticanonical divisor $-K_X$ of $X$ is $\sum_{\gra\in\ol{R}^\cech_{G,\theta}}
a_{\gra}D_\gra+\sum_{E\in N(X)}E$ with $
\sum a_{\gra}\omega_{\gra}=2 \rho-2 \rho_{0}$.
Here  $2\rho:=2\rho_{R_{G}}=\sum_{\alpha\in R_{G}^{+}} \omega_{\alpha}$ is
the sum of all the positive roots of $R_{G}$, while
$2\rho_0:=2\rho_{R_{G}^{0}}$ is the sum of the positive roots in $R_{G}^{0}$.

Let $k$ (or $k_X$) be the piecewise linear function associated to $-K_X$. The anticanonical divisor $-K_X$ is linearly
equivalent to a unique $G$-stable divisor $-\wt{K}_X$. The  piecewise linear function $\widetilde{k}$ (or $\widetilde{k}_X$) associated to $-\wt{K}_X$
is equal to $k-2\rho+2\rho_0$ over $N(X)$ and to 0 over $\rho(D(X))$.  Indeed $-\wt{K}_X$ is $-K_X+div(\prod_{\gra\in \ol{R}_{G,\grt}}f_\gra^{a_\gra})$,
where $f_\gra\in \mC(G/H)^{(B)}$ is an equation of $D_\gra$ (of weight $\gro_\gra$). 
In particular $\widetilde{k}=k-2\rho+2\rho_0$   if $X$ is toroidal.


\section{Standard symmetric varieties}\label{sez standard}

{\em In the following, unless explicitly stated, we always suppose $G$ semisimple
(we will consider the general reductive case mainly in \S\ref{sez rang2+G red}).
Moreover, we  often denote the normalizator $N_G(H)$ by $N(H)$.}
In this section we show that all the standard symmetric varieties are
quasi $\mQ$-Fano varieties; moreover we classify the Fano ones. When   the rank of $G/H$ is one, the standard compactification  $X_0$ of $G/H$ is the unique $G$-equivariant
compactification. In such a case  $X_0$ is an homogeneous projective variety w.r.t.
$Aut^{0}(X)$ by  \cite{A83}; moreover it is wonderful and Fano, because   either it is $\mP^n\times\mP^n$ or it has Picard number one.

First of all, we  reduce ourselves to the indecomposable case. Write
$(G,\theta)$ as a product $\prod(G_{j},\theta)$ of indecomposable
involutions and let  $X_{j}$ be standard compactification of $G_{j}/(G_{j}\cap H)$.

\begin{lem}\label{wonderful+decomposable}
The variety $X$ is (quasi) $\mQ$-Fano if and only if all the $X_{i}$ are
(quasi) $\mQ$-Fano.
\end{lem}

{\em Proof.} The weight $\widetilde{k}_{X}$ is equal to $\sum \widetilde{k}_{X_{i}}$. $\square$.

A standard symmetric variety is always $\mQ$-factorial; in particular,  $K_X$ is a $\mQ$-Cartier divisor.
Moreover, if $X$ is wonderful then also the $X_j$ are wonderful and $X=\prod X_j$ (see \cite{Ru09} Corollary 2.1).
We have the following theorem:

\begin{thm}\label{wonderful+indecomposable}
Let   $X$ be a standard indecomposable  symmetric variety. Let $n$
be the rank of $G$ and let $l$ be the rank of $X$. Then:
\begin{itemize}
\item The anticanonical divisor of $X$ is always a nef and big $\mQ$-divisor.

\item Suppose $X$   wonderful. Then it is not a Fano variety if
and only: i) if the involution induced on $\chi(S)_{\mR}$ is $-id$; ii) $R_{G,\grt}$ is different from $A_{n}$ and $B_{n}$; iii) $H=N_G(G^\grt)$.

\item The standard  indecomposable    varieties whose
anticanonical divisor is not   ample   are  compactifications of the
symmetric spaces in Figure \ref{fig standard indecom}.

\begin{figure}
\begin{changemargin}{0.5cm}{1cm}
\begin{tabular}{|c|c|c|c|}\hline
$G/H$& $\theta$&$n,l$&wonderful\\ \hline\hline

$Spin_{2n+1}/(Spin_{l}\times Spin_{2n+1-l}$)& $BI$&$n\geq l\geq4$&no\\\hline

$Spin_{7}/(Spin_{3}\times Spin_{4})$& $BI$&
$n=l=3$&no\\\hline

$Sp_{2n}/N(GL_{n})$& $CI$&$n=l\geq3$&yes\\\hline

$Spin_{2l}/N(Spin_{l}\times Spin_{l})$&$DI$&$n=l\geq4$ &yes\\\hline

$Spin_{2n}/H$ &$DI$&$n=l\geq6$ &no\\
$(G^{\grt} \subset H\subsetneq N(G^{\grt}))$&&&\\ \hline

$Spin_{2n}/(Spin_{l}\times Spin_{2n-l})$ &$DI$& $n>l\geq4$&no\\ \hline

$Spin_{8}/(Spin_{3}\times Spin_{5})$&$DI$&$n=4,l=3$ &no\\\hline

$E_{6}/N (C_{4})$& $EI$& &yes\\ \hline

$E_{7}/N (A_{7})$ &$EV$&&yes\\ \hline

$E_{7}/A_{7}$ &$EV$&&no\\ \hline

$E_{8}/D_{8}$&$EVIII$&&yes\\ \hline

$F_{4}/(C_{3}\times A_{1})$&$FI$&&yes\\ \hline

$G_{2}/(A_{1}\times A_{1})$&$G$&&yes\\
\hline
\end{tabular}
\end{changemargin}
\caption{Non $\mQ$-Fano standard indecomposable symmetric varieties}
\label{fig standard indecom}
\end{figure}
\end{itemize}
\end{thm}

{\em Proof.}  We have  to determine when  $(\widetilde{k}_{C^-},\alpha_{i})\leq0$
for each $i\in\{1,...l\}$.  We can write $-2\rho+2\rho_0$ as the sum of the
spherical weights $-2\widetilde{\rho} =-2\sum_{\beta_{j}\in
\overline{R}_{G}^{1}}\varpi_{j}$ and $2\widetilde{\rho}_{0}=$ $
-2\sum_{\beta_{j}\in
\overline{R}_{G}^{0}}\varpi_{j}+2\rho_{0} $. Write
$\beta_{j}-\theta(\beta_{j})=\alpha_{i}$, so     $(\widetilde{k}_{C^-},\gra_i )= 2(k_{C^-},\beta_{j} )$ and $(2\widetilde{\rho}_{0},\alpha_{i})=
4(\rho_{0},\beta_{j})\,\leq0$. Thus $(-2\rho+2\rho_0)(\alpha^{\vee}_{i})=-1$ if $\omega_{i}=2\varpi_{j}$ and
$( -2\rho+2\rho_0)(\alpha_{i}^{\vee})\leq-2$ otherwise.  Suppose now  $H$ autonormalizing, i.e. $H=N(G^{\theta})$; the case where $H\subsetneq N(H)$ is very similar. We want to study $k_{C^-}=-\sum_{i=1}^{l} \alpha_{i}$. By the expression of  the Cartan matrix of $R_{G,\theta}$,  $k_{C^-}(\alpha^{\vee}_{i})\leq 1$ for each $i$. Therefore $\widetilde{k}_{C^-}$ is always anti-dominant.
If $\widetilde{k}_{C^-}$ is not regular, then there is (a unique)
$\alpha_{i}=2\beta_{j}\in \ol{R}_{G,\grt}$ such that $k_X (\alpha_{i}^\cech)=1$; in particular $G$ is
simple.
By the classification of the involutions by their Satake
diagrams, $\widetilde{k}_{C^-}$ is not regular if and only if $\theta=-id $ over
$\chi(S)$ and $R_{G,\theta}$ is different by $A_{l}$ and
$B_{l}$.
$\square$

\section{Blows-ups along closed orbits}

In this sections  we want to prove a partial result in arbitrary
rank. We restrict ourselves to the smooth toroidal case. For the toric
varieties of rank 2 one can easily proves  the following property (*):

{\em Let $Z$ be a smooth toric variety of rank 2 and let $Z'$ be a
smooth toric variety birationnaly proper over $Z$. If the
anticanonical bundle of $Z'$ is ample, then also  the anticanonical
bundle of  $Z$ is ample.}

This allows to prove easily that a smooth toric variety proper over
$\mathbb{A}^{2}$ with ample anticanonical bundle is either
$\mathbb{A}^{2}$ or its blow-up in the origin. We would like to use a similar property to classify the
smooth toroidal Fano symmetric varieties. Unfortunately the
previous property (*) is false already in rank three. Indeed, let $\overline{Z}$
be the 3-dimensional toric variety  whose fan $\overline{\Delta}$
has maximal cones $cone(e_{1},e_{2},e_{3})$ and
$cone(e_{1},e_{2},e_{1}+e_{2}-e_{3})$, where $\{e_{1},e_{2},e_{3}\}$ is any
basis of $\chi_{*}(S)$. The function associated to its
anticanonical bundle is linear, so $-K_{Z}$ is nef but non-ample.
Furthermore, the blow-up $\overline{Z} \,'$ of $\overline{Z}$ along
$cone(e_{1},e_{2})$ is  Fano. We can take $\overline{Z}\times \mathbb{A}^{m}$ as  higher dimensional example.

In the previous example, we have considered a blow-up along a
subvariety with strictly positive dimension. Now, we  prove a property similar to (*) considering only  blow-ups along
compact orbits, i.e. $S$-fixed points. In the next section, we prove a much stronger
statement when the rank is three.

\begin{lem} \label{toric rango almeno 3+prop(*)}
Let  $Z$ be a smooth $l$-dimensional toric variety  whose fan  contains
two $l$-dimensional cones $\sigma_{+}$ and $\sigma_{-}$ such that: i)
$\sigma_{+}\cap\sigma_{-}$ has dimension $l-1$ and ii)
$\sigma_{+}\cup\sigma_{-}$ is strictly convex. Assume moreover that
the piecewise linear function associated to the anticanonical bundle
of $Z$ is not strictly convex on  $\sigma_{+}\cup\sigma_{-}$. Then
the anticanonical bundle of any toric variety obtained from $Z$
by a sequence of blow-ups centred in $S$-fixed points is not
ample.
\end{lem}

We can reformulate the first hypothesis  in a more combinatorial
way. Indeed,   we can write
$\sigma_{+}=cone(v_{1},...,v_{l-1},v_{+})$ and
$\sigma_{-}=cone(v_{1},...,$ $v_{l-1},$ $v_-)$ with
$v_{1},...,v_{l-1},v_{+},v_{-}$ primitive and $v_{+}+v_{-}=\sum
a_{i}v_{i}$, where the $a_{i}$ are positive integers, not all zero.

{\em Proof of Lemma \ref{toric rango almeno 3+prop(*)}.} The anticanonical bundle of a variety with satisfies the hypotheses of the lemma is not ample. Thus, it is sufficient to show that the  blow-up $Z'$  of $Z$ centered in any  $\sigma\in\grD(l)$ satisfies again the hypotheses of the lemma. We can suppose
$\sigma=\sigma_{+}$ by symmetry. Then the fan of $Z'$ contains $\sigma_{-}$ and
$\sigma':=cone(v_{1},...,v_{l-1}, v_{+}+\sum v_{i})$.
We have $(v_{+}+\sum v_{i})+v_{-}= \sum (a_{i}+1)v_{i}$, so $Z'$ satisfies the hypotheses w.r.t. $\sigma_{-}$
and $\sigma'$.  $\square$

Now we can classify the toric varieties with ample anticanonical
bundle which are obtained from $\mA^{l}$ by a sequence
of blow-ups centered in $S$-fixed points.

\begin{prop}\label{toric fano almeno 3}
Let  $Z$ be a smooth  toric variety with ample anticanonical bundle which
is obtained from $\mA^{l}$ by a sequence of blow-ups
centered in $S$-fixed points, then it is either $\mA^{l}$ or
the  blow-up of $\mA^{l}$ in the $S$-fixed point.
\end{prop}

{\em Proof.} One can easily see that  the blow-up $Z_{1}$ of
$\mA^{l}$ in the $S$-stable point has ample anticanonical
bundle. The blow-up of $Z_{1}$ in the $S$-fixed point corresponding to
$cone(e_{1},...,\hat{e}_{j},...,e_{l},\sum e_{i})$ satisfies the
hypotheses of the previous lemma w.r.t.
$cone(e_{1},...,\hat{e}_{h},...,e_{l},\sum_{i=1}^{l} e_{i})$ and
$cone(e_{1},...,\hat{e}_{h},...,\hat{e}_{j}, $ $...,$ $e_{l},$
$\sum_{i=1}^{l} e_{i},$ $ 2\sum_{i=1}^{l}e_{i}-e_{j})$, where $h\neq
j$. $\square$

Thus,  a symmetric variety obtained from a wonderful one by a
sequence of blow-ups along closed orbits can be Fano only if it is either
the
wonderful variety or its blow-up along the closed orbit. We have
already considered the wonderful case. Now, we  prove that, when such
blow-up   is Fano, the rank of every indecomposable
factor of $R_{G,\grt}$ is at most 3.

\begin{lem}\label{bound rango fano blow-up almeno 3}
Let $X_{1,...,l}$ be the blow-up of the wonderful compactification of $G/H$ and
suppose that $R_{G,\theta}$ contains an irreducible factor of rank
at least three. If $X_{1,...,l}$ is Fano then it is indecomposable,  has rank 3 and $H\subsetneq
N(G^{\theta})$.
\end{lem}

{\em Proof.} The weights $\{\widetilde{k}_C\}$ associated to $-\widetilde{K}_{X_{1,...,l}}$ are
$\lambda_{i}=-2\rho+2\rho_{0}-(l-2)e_{i}^{*}+\sum_{j\neq
i} e_{j}^{*}$ with $i=1,...,l$.
First, suppose   $G/H$ indecomposable and write $e_i^*=-x_i\gra_i$. If $H\subsetneq
N(G^{\theta})$, then $R_{G,\grt}$ has type $A_1$, $B_2$ or $C_l$.

We consider two cases. First, suppose that there is
$\beta_{h}\in\overline{R}_{G}^{1}$ orthogonal to $R_{G}^{0}$. Write
$\beta_{h}-\theta(\beta_{h})=\alpha_{j}$, so
$0>(\lambda_{j},\alpha_{j}^{\vee})=
(-2\rho+2\rho_{0},\alpha_{j}^{\vee}))+((l-2)x_{j}\alpha_{j}^{\vee},\alpha_{j}^{\vee})-(\sum_{i\neq
j} x_{i}\alpha_{i}^{\vee},\alpha_{j}^{\vee})\geq-2+2(l-2)x_{j}+0$. Observe that   $x_j^{-1}\leq2$, so $l=3$ and  $H\subsetneq N(G^{\theta})$.

If there is not  such a root, we have
the following possibilities:   1) $\theta$ has type  $AII$ and
$G/G^{\theta}$ is $SL_{2l+2}/Sp_{2l+2}$; 2) $\theta$ has  type $CII$ and
$G/G^{\theta}$ is $Sp_{2n}/(Sp_{2l}\times Sp_{2n-2l})$; 3) $\theta$
has type $DIII$ and $G/G^{\theta}$ is $SO_{4l+2}/GL_{2l+1}$. Then there are
$\beta_{3},\beta_{5}\in R_{G}^{0}$ orthogonal to $\overline{R}^0_G\backslash\{\grb_3,\grb_5\}$ and
$ \beta_{4} \in R_{G}^{1}$ such that
$\alpha_{2}=\beta_{3}+2\beta_{4}+\beta_{5}$. Moreover,
$(\beta_{3},\beta_{3})=(\beta_{4},\beta_{4})=(\beta_{5},\beta_{5}) $, $\gra_2^\cech=\frac{1}{(\beta_{2},\beta_{2})}\gra_2$
 and $x_i=1$ if $i<l$.
Thus $0>(\lambda_{2},\alpha_{2}^{\vee})=
(-2\rho+2\rho_{0},\alpha_{2}^{\vee}))+(l-2)(\alpha_{2},\alpha_{2}^{\vee})
-(\alpha_{1} ,\alpha_{2}^{\vee})-x_{3}(\alpha_{3}
,\alpha_{2}^{\vee})\geq-4+2(l-2)+1+x_{3}$, so  we have again
$l=3$  and $H\subsetneq N(G^{\theta})$.

Finally, suppose   $\theta$   decomposable.  Let
$(G,\theta)=(G_{1},\theta_{1})\times(G_{2},\theta_{2})$ with $l':=rank\,
(G_{1}/G_{1}^{\theta})\geq3$ and define the weight $\lambda_{i}'$ for
$G_{1}$ in an analogous way to the $\lambda_{i}$. We have
$\lambda_{i}=\lambda'_{i}-(l-l')e_{i}^{*}+\omega$
where $\omega$ is   orthogonal to
$R_{G_{1},\theta}$. By the previous part of the proof  there is
always an $i$ with  $\lambda_{i}'(\alpha_{i}^{\vee})\geq-1$, so
$\lambda_{i}(\alpha_{i}^{\vee})\geq \lambda_{i}'(\alpha_{i}^{\vee})+ \frac{1}{2}
(\alpha_{i},\alpha_{i}^{\vee})\geq0$, a contradiction.   $\square$

By an explicit analysis of   the indecomposable involutions of rank at most three we obtain:

\begin{thm}\label{fano rango almeno 3} Let $G/H$ be a symmetric space of rank $l$ ($>1$) associated to an involution $\grt$ and let $X$ be a compact symmetric variety  obtained  from the wonderful compactification of $G/H$  by a sequence
of blow-ups along closed orbits.
\begin{enumerate}
\item If $X$ is a Fano variety then either it is the wonderful variety $X_{0}$ or it is the blow-up $X_{1,...,l}$ of $X$ along the
closed orbit.
\item If there is an indecomposable   factor of $(G,\theta)$ of rank at least 3 then $X_{1,...,l}$ is not Fano.
\item If $(G,\theta)$ has rank at least 6 and has an indecomposable factor of rank
2, then $X_{1,...,l}$ is not Fano.
\item If $X_{1,...,l}$ is Fano, the possibilities
for the indecomposable factors of $G/H$ are as in Figure \ref{figura fano at least 3} (we indicate also the eventual
conditions on $rank\,G/H$ so that such a factor  can appear).
\end{enumerate}
\end{thm}


\section{Regular Fano varieties of rank 3}\label{sez rk3}

In this section  we suppose $X_0$ wonderful; recall that $\{e_1,e_2,e_3\}$ is the
basis of $\chi_*(S)$ which generates  $C^-$. We classify  all the Fano symmetric varieties
obtainable by a wonderful symmetric variety of rank three from $X_0$ by a
succession of blow-ups along $G$-subvarieties. This class of
varieties contains many varieties; indeed each compact symmetric
variety is dominated by a smooth toroidal variety obtained by a
succession of blow-ups along $G$-subvarieties of codimension two. We
begin proving a result similar to   Lemma \ref{toric rango almeno
3+prop(*)}.

\begin{figure}\begin{changemargin}{0.5cm}{1cm}
\begin{tabular}{|c|c|c|c|}\hline
$G'/H'$&$\theta|G'$& $rank\,G/H$ &$rank\,G'/H'$ \\
\hline\hline

$SL_{3} $&$A_2$& $l=2$&2\\\hline

$Spin_{5}$&$B_2$& $l=2$&2\\\hline

$SL_{6}/N(Sp_{6})$&$AII$& $l=3$&2\\\hline

$SL_{n+1}/S(GL_2\times GL_{n-1}) $, $n\geq4$&$AIII$& $l=2$&2\\\hline

$Sp_{2n}/Sp_4\times Sp_{2n-4}$, $n\geq5$&$CII$&$l\leq3$&2\\\hline

$Sp_8/Sp_4\times Sp_{4}$&$CII$& $l=3$&2\\\hline

$Sp_{8}/N(Sp_{4}\times Sp_{4})$&
$CII$&$l=2$&2\\\hline

$SO_{10}/GL_{5}$& $DIII$& $l\leq3$&2\\\hline

$E_{6}/N(D_{5}\times \mathbb{C}^{*})$& $EIII$& $l\leq4$&2\\\hline

$E_{6}/N(F_{4})$& $EIV$& $l\leq5$&2\\\hline

$PSL_{2}$&$A_{1}$& $l=2$&1\\\hline

$SL_{2}$&$A_{1}$&$l\leq3$&1\\\hline

$SL_{n+1}/N(SO_{n})$& $AI$ & $l=2$&1\\\hline

$SL_{n+1}/SO_{n}$& $AI$ & $l=2$&1\\\hline

$SL_{6}/N(Sp_{6})$& $AII$& $l\leq3$&1\\\hline

$SL_{6}/Sp_{6}$& $AII$& $l\leq5$&1\\\hline

$SL_{n+1}/N(S(GL_{1}\times GL_{n}))$& $AIV$ & $l\leq n+1$&1\\\hline
$SO_{2n+1}/N(SO_{1}\times SO_{2n})$ & $BII$ &$l\leq n+1$&1\\\hline
$SO_{2n+1}/(SO_{1}\times SO_{2n})$ & $BII$ &$l\leq 2n$&1\\\hline
$Sp_{2n}/(Sp_{2}\times Sp_{2n-2})$& $CII$ & $l\leq 2n$&1\\\hline
$SO_{2n}/N(SO_{1}\times SO_{2n-1})$& $DII$ & $l\leq n$&1\\\hline
$SO_{2n}/(SO_{1}\times SO_{2n-1})$& $DII$ & $l\leq 2n-1$&1\\\hline

$F_{4}/B_{4}$& $FII$ & $l\leq12$&1\\\hline
\end{tabular}\end{changemargin}
\caption{Fano $X_{1,...,l}$} \label{figura fano at least 3}
\end{figure}

Let $\widetilde{Z}$ be the toric variety  whose fan
$\widetilde{\Delta}$ has maximal cones $cone(v_{1},v_{2},v_{+})$ and
$cone(v_{1},v_{2},v_{-}=x_{1}v_{1}+x_{2}v_{2}-v_{+})$, where
$\{v_{1},v_{2},v_{+}\}$ is a basis of $\chi_*(S)$,
$x_{1}+x_{2}>0$ and $x_{1}\geq x_{2}\geq0$. The anticanonical bundle
of $\widetilde{Z}$ is ample if and only if $x_{1}=x_{1}+x_{2}=1$. In
this case, $\widetilde{Z}$ is the blow up of $\mA^{3}$ along
a stable subvariety of codimension 2. Moreover, the anticanonical
bundle of $\widetilde{Z}$ is nef, but non-ample if and only if
$x_{1}+x_{2}=2$. We have two possibilities: either
$v_{+}+v_{-}=v_{1}+v_{2}$ or $v_{+}+v_{-}=2v_{1}$. In the first case
we have a variety isomorphic to the variety $\overline{Z}$ of the
previous section. This is the more problematic case, so we will study it in a second time.

\begin{lem}\label{toric fano rango 3, lem1}
Let  $Z$ be a smooth 3-dimensional toric variety whose fan contains two maximal
cones $cone(v_{1},v_{2},v_{+})$ and $cone(v_{1},v_{2},v_{-})$ such that
$v_{+}+v_{-}=x_{1}v_{1}+x_{2}v_{2}$, where $x_{1}$ and $x_{2}$ are
integers with $x_{1}\geq x_{2}\geq0$. Suppose moreover that
$x_{1}\geq2$. Then the anticanonical bundle of any toric variety
$Z'$ obtained from $Z$ by a sequence of blow-ups along $S$-subvarieties is not ample.
\end{lem}

{\em Proof.} Remark that the anticanonical bundle of $Z$ is not ample.  We say that a variety satisfies weakly the hypotheses of
the lemma   if $x_{1}+x_{2}\geq2$ (instead of $x_1\geq2$). We use
the following trivial observation: $x_{1}+x_{2}>2$ implies
$x_{1}\geq2$.
One can try to prove this lemma by induction as the
Lemma~\ref{toric fano almeno 3}. Unfortunately we can prove only the following weaker statement.

\begin{lem}\label{lem sub Zbarra}
Let $Z$ be a toric variety which satisfies weakly the hypotheses of Lemma \ref{toric fano rango 3, lem1} and  let $Z'$ be the blow-up of $Z$ along a cone $\tau$.
\begin{enumerate}
\item If $\tau\neq cone(v_{1},v_{2})$, then $Z'$ satisfies   the hypotheses
of Lemma \ref{toric fano rango 3, lem1}.
\item If $Z$ satisfies  the hypotheses of Lemma \ref{toric fano rango 3, lem1}, then $Z'$ satisfies
weakly the hypotheses of Lemma \ref{toric fano rango 3, lem1}.
\end{enumerate}
\end{lem}

{\em Proof.} We can suppose $\tau\subset cone(v_{1},v_{2},v_{-})$ by symmetry. If $\tau\neq
 cone(v_{1},$ $v_{2})$, we have three
possibilities: $\tau=cone(v_{1},v_{-})$, $\tau=cone(v_{2},v_{-})$
and $\tau=cone(v_{1},v_{2},v_{-})$. We  always have
$\Delta_{Z'}[p]=\Delta_Z[p]\cup\{v':=v_{-}+b_{1}v_{1}+b_{2}v_{2}\}$
with $b_{1},b_{2}\in\{0,1\}$. Moreover, $\Delta_{Z'}$ contains the cones
$cone(v_{1},v_{2},v_{+})$ and
$cone(v_{1},v_{2},v')$ and we have
$v'+v_{+}=(x_{1}+b_{1})v_{1}+(x_{2}+b_{2})v_{1}$
with $(x_{1}+b_{1})+(x_{2}+b_{2})>2$, so $Z_{1}$ satisfies the
hypotheses of the lemma.

Finally let  $\tau=cone(v_{1},v_{2})$. The fan of $Z'$ contains the
cones $cone(v_{1},v_{1}+v_{2},v_{+})$ and
$cone(v_{1},v_{1}+v_{2},v_{-})$. We have
$v_{+}+v_{-}=(x_{1}-x_{2})v_{1}+x_{2}(v_{1}+v_{2})$ with
$(x_{1}-x_{2})+x_{2}=x_{1}\geq2$.   $\square$

\

Now, we consider the general case. We have a sequence $Z=Z_{0}\leftarrow
Z_{1}\leftarrow...\leftarrow Z_{i}\leftarrow ...\leftarrow Z_{r}=Z'$
where $Z_{i+1}$ is the blow-up of $Z_{i}$ along the cone $\tau_{i}$.  Let $\Delta_{i}=\grD_{Z_i}$  and let $j$ be the maximal index such that  $Z_{j}$
satisfies the hypotheses (w.r.t.
$cone(w_{1},w_{2},w_{+})$ and
$cone(w_{1},w_{2},w_{-})$). By the previous lemma $Z_{j+1}$   satisfies weakly the
hypotheses, in particular its anticanonical bundle is not ample. By the maximality of $j$,  $Z_{j+1}$ does not satisfies the hypotheses, so $Z_{j+1}$ contains a
variety isomorphic to $\overline{Z}$. Let $\overline{ \Delta}$ be
the fan of such variety. If $\grD_r$ contains $\overline{\grD}$ then $-K_{Z'}$  is not ample. Otherwise there is a minimal $h$
such that $\overline{ \Delta}$ is not contained in $\grD_{h+1}$. We claim that   $Z_{h+1}$ satisfies the
hypotheses of the lemma, a contradiction.

By the previous lemma  $\tau_{j}=cone(w_{1},w_{2})$.
We know that $Z_{j+1}$ satisfies
weakly the hypotheses w.r.t.
$\sigma_{+}=cone(w_{1},w_{1}+w_{2},w_{+})$ and
$\sigma_{-}=cone(w_{1},w_{1}+w_{2},w_{-})$. Moreover, the open
subvariety of $Z_{j+1}$ with maximal cones $\sigma_{+}$ and
$\sigma_{-}$  is isomorphic to $\overline{Z}$.
By the Lemma \ref{lem sub Zbarra}  2) we can suppose  $\tau_{h}=cone(w_{1},w_{1}+w_{2})$.

The   fan of   $Z_{h+1}$   contains two cones
$cone(w_{1}+w_{2},w_{+},2w_{1}+w_{2})$ and $
cone(w_{1}+w_{2},$ $w_{+},$ $w')$, with $w'=w_2+x_1(w_1+w_2)+x_2w_+$ and $x_1,x_2\geq0$. Therefore
$Z_{h+1}$ satisfies the hypotheses of the lemma w.r.t.
these cones. Indeed  $(2w_{1}+w_{2})+w'=
(2+x_1)(w_{1}+w_{2})+x_2w_+$. $\square$

Now we want to study the varieties which contain an open subvariety
isomorphic to $\overline{Z}$. Observe that these varieties are never
Fano varieties. Let  $Z$ be  such a variety and let $Z'$ be the
blow-up   of $Z$ along  the subvariety
of $\overline{Z}$ associated to $cone(e_{1},e_{2})$. We prove
that, if $Z'$ satisfies the hypotheses of  Lemma \ref{toric fano rango 3, lem1},
then there are not Fano varieties obtainable from $Z$ by a
sequence of blow-ups.

\begin{lem}\label{toric fano rango 3, lem2}
Let  $Z$ be a smooth 3-dimensional toric variety whose fan contains  $cone(v_{1},$ $v_{2},v_{3})$ and
$cone(v_{1},v_{2},v_{1}+v_{2}-v_{3})$ for suitable
$v_{1},v_{2},v_{3}$. Let  $Z'$ be the blow-up of
$Z$ along the stable subvariety corresponding to $cone(v_{1},v_{2})$
and let  $Z''$ be a toric variety obtained from $Z$ by a
sequence of blow-ups along $S$-subvarieties. If the anticanonical bundle of $Z''$ is ample,
then $Z''$ is obtainable from $Z'$ by a sequence of blow-ups along $S$-subvarieties.
\end{lem}

{\em Proof.} We cannot proceed as in the previous lemma, because we
do not know  the other  cones  of $\grD_Z$. We have again a sequence
$Z=Z_{0}\leftarrow Z_{1}\leftarrow...\leftarrow Z_{i}\leftarrow
...\leftarrow Z_{h}=Z''$ where $\pi_{i+1}:Z_{i+1}\rightarrow Z_{i}$
is the blow-up  along  $\tau_{i}$.  First of all, there is a (minimal) cone $\tau_{j}$
contained in $cone(v_{1},v_{2},v_{3},v_{1}+v_{2}-v_{3})$, because
otherwise the anticanonical bundle of  $Z''$ is not ample. By Lemma  \ref{toric
fano rango 3, lem1}, $\tau_{j}$ is $cone(v_{1},v_{2})$. We want to reorder the cones associated to the
subvarieties along which we are blowing-up. Clearly this
operation is not well defined  in general.

We  consider the following sequence of blow-ups:
$Z=Z'_{0}\leftarrow Z'_{1}\leftarrow...\leftarrow Z_{i}'\leftarrow
...\leftarrow Z'_{j+1}$, where $\pi'_{0}:Z'_{1}\rightarrow Z'_{0}$  is
the blow-up  along
$\tau_{j}$  and $\pi'_{i+1}:Z'_{i+1}\rightarrow Z'_{i}$ is the blow-up along   $\tau_{i-1}$ for
each $i\geq1$. Let   $\Delta'_{i}=\grD_{ Z'_{i}}$. We show that these blow-ups are well defined and that
$Z'_{j+1}=Z_{j+1}$.

The cone $\tau_{i-1}$ belongs to $\Delta'_{i}$ for each
$1\leq i\leq j$  because $\tau_{i}$ is
contained in $|\Delta| \backslash
cone(v_{1},v_{2},v_{3},v_{1}+v_{2}-v_{3})$ for each $i\leq j$. Moreover,
the elements of $\Delta'_{1}(3)$ not contained in
$cone(v_{1},v_{2},$ $v_{3},$ $v_{1}+v_{2}-v_{3})$ are exactly the
elements of $\Delta_{0}(3)\setminus\{cone(v_{1},v_{2},v_{3}),
cone(v_{1},v_{2},v_{1}+v_{2}-v_{3})\}$.

$Z$ is the union of the following two open $S$-subvarieties: $U_{1}$
whose  fan  has maximal cones $cone(v_{1},v_{2},v_{3})$ and
$cone(v_{1},v_{2},v_{1}+v_{2}-v_{3})$; $U_{2}$ whose fan has maximal
cones $\Delta(3)\setminus\{cone(v_{1},v_{2},v_{3}),
cone(v_{1},v_{2},v_{1}+v_{2}-v_{3})\}$.

The blow-up $\pi'_{0}$ induces an isomorphism between $U_{2}$ and
its inverse image, because
$cone(v_{1},v_{2})$ is not contained  in any maximal cone of
$U_{2}$. In the same way $\pi_{j}$ induces an isomorphism between
the inverse image of $U_{2}$ in $Z_{j}$ and its inverse image in
$Z_{j+1}$. So the inverse image of $U_{2}$ in $Z_{j+1}$ is isomorphic to
the the inverse image of  $U_{2}$ in $Z'_{j+1}$. Moreover
$\pi'_{j}\circ...\circ\pi'_{2}$ induces an isomorphism between
$(\pi'_{1})^{-1}(U_{1})$ and its inverse image. In the same way
$\pi_{j-1}\circ...\circ\pi_{1}$ induces an isomorphism between
$U_{1}$ and its inverse image. So the inverse image of $U_{1}$ in
$Z_{j+1}$ is isomorphic to the the inverse image of $U_{1}$ in
$Z'_{j+1}$. The lemmas follows because there is at most one morphism between
two toric $S$-varieties extending  the identity automorphism of $S$. $\square$

We now restrict the   possible Fano   symmetric varieties with rank three (and fixed $G/H$) which are obtainable as before
to a finite explicit list.

\begin{prop}\label{fano toric rango 3}
The   toric varieties obtainable from $\mA^{3}$ by a
sequence of blow-ups and with ample anticanonical bundle are, up to isomorphisms:

\begin{enumerate}

\item $\mA^{3}$;

\item a  2-blow-up of $\mA^{3}$;

\item the 3-blow-up of $\mA^{3}$;

\item the variety whose fan has maximal cones:  $cone(e_{1},e_{1}+e_{2},e_{1}+e_{2}+e_{3})$,
$cone(e_{1},e_{3},e_{1}+e_{2}+e_{3})$,
$cone(e_{2},e_{3},e_{1}+e_{2}+e_{3})$ and
$cone(e_{2},e_{1}+e_{2},e_{1}+e_{2}+e_{3})$. This variety is
obtainable from $\mA^{3}$ by two consecutive blow-ups
along subvarieties of codimension two;

\item the variety whose fan has maximal cones:
$cone(e_{1},e_{3},e_{1}+e_{2}+2e_{3})$,
$cone(e_{1},$ $e_{1}+e_{2}+e_{3},e_{1}+e_{2}+2e_{3})$,
$cone(e_{1},e_{2},e_{1}+e_{2}+e_{3})$,
$cone(e_{2},e_{1}+e_{2}+e_{3},e_{1}+e_{2}+2e_{3})$ and
$cone(e_{2},e_{3},e_{1}+e_{2}+2e_{3})$. This variety is obtainable
from $\mA^{3}$ by a 3-blow up followed by a 2-blow up.

\end{enumerate}

\end{prop}

{\em Proof.} We proceed as follows: the   anticanonical bundle of
$\mA^{3}$ is ample, so we consider all the possible blow-ups
of $\mA^{3}$. Let $Z$ be a blow-up of $\mA^{3}$: 1) if
$Z$ satisfies the hypotheses of Lemma~\ref{toric fano rango 3, lem1}
we know that there are not toric  varieties with ample anticanonical
bundle and obtainable from $Z$ by a sequence of blow-ups; 2) if
$Z$ satisfies the hypotheses of the Lemma~\ref{toric fano rango 3,
lem2} we study the variety $Z'$ of such lemma; 3) finally, if
the anticanonical bundle $Z$ is ample, we reiterate the procedure.
Observe that a priori it is  possible that $Z$ belongs to none  of
the previous cases.
In the following, if two blow-ups of a given variety are isomorphic,
we examine only one of them. Given a toric variety  $Z$,  let $\Delta$ be its  fan and let $k$ be
the piecewise linear function associated to   $-K_Z$. Suppose that
all the maximal cones in $\Delta$ are 3-dimensional. Remember that $-K_Z$ is ample  if and only if, given any cone $C\in\Delta(3)$ and
any   $v\in\Delta[p]$ with $v\notin  C$,
$(k_C )(v)<1$.

Let $Z_{0}=\mA^{3}$; it has ample canonical
bundle. Up to isomorphisms there are two blow-ups of
$\mA^{3}$:  the blow-up $Z_{1}$  along $cone(e_{1},e_{2})$
and the blow-up $Z_{2}$ along $cone(e_{1},e_{2},e_{3})$.

One can show that the anticanonical bundle of $Z_{1}$ is ample
because $\grD_1(2)=\{cone(e_{1},e_{3},e_{1}+e_{2}),
cone(e_{2},e_{3},e_{1}+e_{2})\}$ and $\Delta[p]=\{e_{1},e_{2},e_{3},e_{1}+e_{2}\}$.
The blow-ups of $Z_{1}$ are, up to isomorphisms: i) the blow-up
$Z_{11}$ along $cone(e_{1},e_{3})$; ii)   the blow-up $Z_{12}$ along
$cone(e_{1},e_{1}+e_{2})$; iii) the blow-up $Z_{13}$  along
$cone(e_{1}+e_{2},e_{2},e_{3})$;  iv) the blow-up $Z_{14}$  along
$cone(e_{1}+e_{2},e_{3})$.

The variety $Z_{11}$ satisfies the hypotheses of the
Lemma~\ref{toric fano rango 3, lem2} w.r.t.
$cone(e_{3},e_{1}+e_{2},e_{1}+e_{3})$ and
$cone(e_{3},e_{1}+e_{2},e_{2})$. Hence we have to study the
blow-up $Z_{11b}$ of $Z_{11}$ along $cone(e_{3},e_{1}+e_{2})$. This
variety satisfies the hypotheses of the Lemma~\ref{toric fano rango
3, lem2} w.r.t. $cone(e_{1}+e_{2},e_{1}+e_{3},e_{1})$ and
$cone(e_{1}+e_{2},e_{1}+e_{3},e_{1}+e_{2}+e_{3})$. Hence we have to
study the the blow-up $Z_{11c}$ of $Z_{11b}$  along
$cone(e_{1}+e_{2},e_{1}+e_{3})$. $Z_{11c}$ satisfies the hypotheses
of the Lemma~\ref{toric fano rango 3, lem2} w.r.t.
$cone(e_{1}+e_{2},e_{1}+e_{2}+e_{3},e_{2})$ and
$cone(e_{1}+e_{2},e_{1}+e_{2}+e_{3},2e_{1}+e_{2}+e_{3})$. Hence we
have to consider the blow-up $Z_{11d}$ of $Z_{11c}$  along
$cone(e_{1}+e_{2},e_{1}+e_{2}+e_{3})$. $Z_{11d}$ satisfies the
hypotheses of the Lemma~\ref{toric fano rango 3, lem1} w.r.t.
$cone(e_{2},e_{1}+e_{2}+e_{3},e_{3})$ and
$cone(e_{2},e_{1}+e_{2}+e_{3},2e_{1}+2e_{2}+e_{3})$. Thus there are
not toric  varieties with ample anticanonical bundle and obtained
from $Z_{11}$ by a sequence of blow-ups.

$Z_{12}$  satisfies the hypotheses of the Lemma~\ref{toric fano
rango 3, lem1} w.r.t. $cone(e_{3},e_{1}+e_{2},2e_{1}+e_{2})$ and
$cone(e_{3},e_{1}+e_{2},e_{2})$.  $Z_{13}$ satisfies the hypotheses
of the Lemma~\ref{toric fano rango 3, lem1} w.r.t.
$cone(e_{3},e_{1}+e_{2},e_{1})$ and
$cone(e_{3},e_{1}+e_{2},e_{1}+2e_{2}+e_{3})$.

The anticanonical bundle of $Z_{14}$ is ample because $\grD_{14}(2)=\{
cone(e_{1}+e_{2},e_{1},e_{1}+e_{2}+e_{3}),
cone(e_{1},e_{3},e_{1}+e_{2}+e_{3}),
cone(e_{2},e_{3},e_{1}+e_{2}+e_{3}),
cone(e_{1}+e_{2},e_{2},e_{1}+e_{2}+e_{3})\} $ and $\Delta[p]=\{e_{1},e_{2},e_{3},e_{1}+e_{2},e_{1}+e_{2}+e_{3}\}$.
The blow-ups of $Z_{14}$ are:  i) the blow-up $Z_{141}$ of $Z_{14}$
along  $cone(e_{1},e_{1}+e_{2})$; ii) the blow-up $Z_{142}$ of
$Z_{14}$ along $cone(e_{1}+e_{2},e_{1}+e_{2}+e_{3})$; iii)  the
blow-up  $Z_{143}$ of  $Z_{14}$ along
$cone(e_{3},e_{1}+e_{2}+e_{3})$; iv) the blow-up $Z_{144}$ of
$Z_{14}$ along $cone(e_{2},e_{3},e_{1}+e_{2}+e_{3})$; v)   the
blow-up $Z_{145}$ of $Z_{14}$  along
$cone(e_{2},e_{1}+e_{2},e_{1}+e_{2}+e_{3})$; vi) the blow-up
$Z_{146}$ of  $Z_{14}$  along $cone(e_{1},e_{3})$; vii) the blow-up
$Z_{147}$ of  $Z_{14}$  along $cone(e_{1},e_{1}+e_{2}+e_{3})$.

The variety $Z_{141}$  satisfies the hypotheses of the
Lemma~\ref{toric fano rango 3, lem1} w.r.t.
$cone(e_{1}+e_{2},e_{1}+e_{2}+e_{3},e_{2})$ and
$cone(e_{1}+e_{2},e_{1}+e_{2}+e_{3},2e_{1}+e_{2})$. The variety
$Z_{142}$ satisfies the hypotheses of the Lemma~\ref{toric fano
rango 3, lem1} w.r.t.
$cone(e_{1},e_{1}+e_{2}+e_{3},2e_{1}+2e_{2}+e_{3})$ and
$cone(e_{1},e_{1}+e_{2}+e_{3},e_{3})$. The variety $Z_{143}$
satisfies the hypotheses of the Lemma~\ref{toric fano rango 3, lem1}
w.r.t.  $cone(e_{1},e_{1}+e_{2}+e_{3},e_{1}+e_{2})$ and
$cone(e_{1},e_{1}+e_{2}+e_{3},e_{1}+e_{2}+2e_{3})$. The variety
$Z_{144}$ satisfies the hypotheses of the Lemma~\ref{toric fano
rango 3, lem1} w.r.t. $cone(e_{3},e_{1}+e_{2}+e_{3},e_{1})$
and $cone(e_{3},e_{1}+e_{2}+e_{3},e_{1}+2e_{2}+2e_{3})$. The variety
$Z_{145}$ satisfies the hypotheses of the Lemma~\ref{toric fano
rango 3, lem1} w.r.t. $cone(e_{2},e_{1}+e_{2}+e_{3},e_{3})$
and $cone(e_{2},e_{1}+e_{2}+e_{3},2e_{1}+3e_{2}+e_{3})$.

The variety $Z_{146}$ satisfies the hypotheses of the
Lemma~\ref{toric fano rango 3, lem2} w.r.t.
$cone(e_{1},e_{1}+e_{2}+e_{3},e_{1}+e_{3})$ and
$cone(e_{1},e_{1}+e_{2}+e_{3},e_{1}+e_{2})$. Hence we have to study
the blow-up of $Z_{146}$ along $cone(e_{1},e_{1}+e_{2}+e_{3})$. This
variety is $Z_{11c}$, so there are no toric varieties with ample
anticanonical bundle which are obtained from $Z_{146}$ by a sequence
of blow-ups.

The variety $Z_{147}$ satisfies the hypotheses of the
Lemma~\ref{toric fano rango 3, lem2} w.r.t.
$cone(e_{1}+e_{2},e_{1}+e_{2}+e_{3},e_{2})$ and
$cone(e_{1}+e_{2},e_{1}+e_{2}+e_{3},2e_{1}+e_{2}+e_{3})$. Hence we
have to study the  blow-up $Z_{147b}$ of $Z_{147}$  along
$cone(e_{1}+e_{2},e_{1}+e_{2}+e_{3})$. This variety satisfies the
hypotheses of the Lemma~\ref{toric fano rango 3, lem1}  w.r.t. $cone(e_{2},e_{1}+e_{2}+e_{3},e_{3})$ and
$cone(e_{2},e_{1}+e_{2}+e_{3},2e_{1}+2e_{2}+e_{3})$. Observe that we
have classified the toric  varieties with anticanonical bundle which are
obtained from $Z_{1}$ by a sequence of blow-ups.

The anticanonical bundle of $Z_{2}$ is ample because $\grD_2(2)=\{cone(e_{2},e_{3},e_{1}+e_{2}+e_{3}),
cone(e_{1},e_{3},e_{1}+e_{2}+e_{3}),
cone(e_{1},e_{2},e_{1}+e_{2}+e_{3})\}$ and $\Delta_{2}( p)= \{
e_{1},e_{2}, e_{3},$ $e_{1}+e_{2}+e_{3}\}$. The blow-ups of $Z_{2}$ are, up to isomorphisms: i) the variety
$Z_{14}$; ii) the blow-up $Z_{21}$ of $Z_{2}$ along
$cone(e_{1},e_{2},e_{1}+e_{2}+e_{3})$; iii) the blow-up $Z_{22}$ of
$Z_{2}$ along $cone(e_{3},e_{1}+e_{2}+e_{3})$.

$Z_{21}$ satisfies the hypotheses of the Lemma~\ref{toric fano rango
3, lem1}  w.r.t. $cone(e_{2},e_{1}+e_{2}+e_{3},e_{3})$ and
$cone(e_{2},e_{1}+e_{2}+e_{3},2e_{1}+2e_{2}+e_{3})$.

The anticanonical bundle of $Z_{22}$ is ample because $\grD_{22}(2)=\{cone(e_{1},e_{3},e_{1}+e_{2}+2e_{3}),
cone(e_{1},e_{1}+e_{2}+e_{3},e_{1}+e_{2}+2e_{3}),
cone(e_{1},e_{2},e_{1}+e_{2}+e_{3}),
cone(e_{2},e_{1}+e_{2}+e_{3},e_{1}+e_{2}+2e_{3}),
cone(e_{2},e_{3},e_{1}+e_{2}+2e_{3}) \}$ and $\Delta_{22}[p]= \{
e_{1},e_{2},e_{3},e_{1}+e_{2}+e_{3},e_{1}+e_{2}+2e_{3}\} $.

The blow-ups of $Z_{22}$  are, up to isomorphisms: i) the variety
$Z_{143}$ which satisfies the hypotheses of the Lemma~\ref{toric
fano rango 3, lem1}; ii)  the blow-up $Z_{221}$ of $Z_{22}$ along
$cone(e_{3},e_{1}+e_{2}+2e_{3})$; iii)  the blow-up $Z_{222}$ of
$Z_{22}$  along $cone(e_{1}+e_{2}+e_{3},e_{1}+e_{2}+2e_{3})$; iv)
the blow-up $Z_{223}$ of $Z_{22}$  along
$cone(e_{1},e_{1}+e_{2}+2e_{3})$; v)  the blow-up $Z_{224}$ of
$Z_{22}$  along $cone(e_{1},e_{3})$; vi)   the blow-up $Z_{225}$ of
$Z_{22}$  along $cone(e_{1},e_{1}+e_{2}+e_{3})$; vii)  the blow-up
$Z_{226}$ of  $Z_{22}$ along $cone(e_{1},e_{3},e_{1}+e_{2}+2e_{3})$;
viii)  the blow-up $Z_{227}$ of $Z_{22}$  along
$cone(e_{1},e_{1}+e_{2}+e_{3},e_{1}+e_{2}+2e_{3})$; ix)   the
blow-up $Z_{228}$ of $Z_{22}$ along
$cone(e_{1},e_{2},e_{1}+e_{2}+e_{3})$.

$Z_{221}$ satisfies the hypotheses of the Lemma~\ref{toric fano
rango 3, lem1} w.r.t.
$cone(e_{1},e_{1}+e_{2}+2e_{3},e_{1}+e_{2}+e_{3})$ and
$cone(e_{1},e_{1}+e_{2}+2e_{3},e_{1}+e_{2}+3e_{3})$. $Z_{222}$
satisfies the hypotheses of the Lemma~\ref{toric fano rango 3, lem1}
w.r.t. $cone(e_{1},e_{1}+e_{2}+2e_{3},e_{3})$ and
$cone(e_{1},e_{1}+e_{2}+2e_{3},2e_{1}+2e_{2}+3e_{3})$. $Z_{223}$
satisfies the hypotheses of the Lemma~\ref{toric fano rango 3, lem1}
w.r.t. $cone(e_{1}+e_{2}+e_{3},e_{1}+e_{2}+2e_{3},e_{2})$
and $cone(e_{1}+e_{2}+e_{3},e_{1}+e_{2}+2e_{3},2e_{1}+e_{2}+2e_{3})$.

$Z_{224}$ satisfies the hypotheses of the Lemma~\ref{toric fano
rango 3, lem2} w.r.t.
$cone(e_{1},e_{1}+e_{2}+2e_{3},e_{1}+e_{2}+e_{3})$ and
$cone(e_{1},e_{1}+e_{2}+2e_{3},e_{1}+e_{3})$. Hence we have to study
the blow-up $Z_{224b}$ of $Z_{224}$  along
$cone(e_{1},e_{1}+e_{2}+2e_{3})$. This variety satisfies the
hypotheses of the Lemma~\ref{toric fano rango 3, lem1} w.r.t. $cone(e_{1},e_{1}+e_{2}+e_{3},e_{2})$ and
$cone(e_{1},e_{1}+e_{2}+e_{3},2e_{1}+e_{2}+2e_{3})$. $Z_{225}$
satisfies the hypotheses of the Lemma~\ref{toric fano rango 3, lem2}
w.r.t. $cone(e_{1},e_{1}+e_{2}+2e_{3},2e_{1}+e_{2}+e_{3})$
and $cone(e_{1},e_{1}+e_{2}+2e_{3},e_{3})$. Hence we have to study
the blow-up $Z_{225b}$  of $Z_{225}$ along
$cone(e_{1},e_{1}+e_{2}+2e_{3})$. $Z_{225b}$  satisfies the
hypotheses of the Lemma~\ref{toric fano rango 3, lem2}  w.r.t.  $cone(2e_{1}+e_{2}+e_{3},e_{1}+e_{2}+2e_{3},e_{1}+e_{2}+e_{3})$
and $cone(2e_{1}+e_{2}+e_{3},e_{1}+e_{2}+2e_{3},2e_{1}+e_{2}+2e_{3})$.
Hence we have to consider the blow-up $Z_{225c}$ of $Z_{225b}$ along
$cone(2e_{1}+e_{2}+e_{3},e_{1}+e_{2}+2e_{3})$. This variety
satisfies the hypotheses of the Lemma~\ref{toric fano rango 3, lem1}
w.r.t.  $cone(e_{1}+e_{2}+e_{3},e_{1}+e_{2}+2e_{3},e_{2})$
and $cone(e_{1}+e_{2}+e_{3},e_{1}+e_{2}+2e_{3},3e_{1}+2e_{2}+3e_{3})$.

$Z_{226}$ satisfies the hypotheses of the Lemma~\ref{toric fano
rango 3, lem1} w.r.t.
$cone(e_{1},e_{1}+e_{2}+2e_{3},e_{1}+e_{2}+e_{3})$ and
$cone(e_{1},e_{1}+e_{2}+2e_{3},2e_{1}+e_{2}+3e_{3})$. $Z_{227}$
satisfies the hypotheses of the Lemma~\ref{toric fano rango 3, lem1}
w.r.t. $cone(e_{1}+e_{2}+e_{3},e_{1}+e_{2}+2e_{3},e_{2})$
and $cone(e_{1}+e_{2}+e_{3},e_{1}+e_{2}+2e_{3},3e_{1}+2e_{2}+3e_{3})$.
$Z_{228}$ satisfies the hypotheses of the Lemma~\ref{toric fano
rango 3, lem1} w.r.t.
$cone(e_{1},e_{1}+e_{2}+e_{3},e_{1}+e_{2}+2e_{3})$ and
$cone(e_{1},e_{1}+e_{2}+e_{3},2e_{1}+2e_{2}+e_{3})$. $\square$

Let $G/H$ be a symmetric space of rank three such that the standard
compactification $X_{0}$ of $G/H$ is wonderful. We introduce the following
notations:
\begin{itemize}
\item we denote by $X_{ij}$ the blow-up of $X_{0}$ along
$(cone( e_{i}, e_{j}),\vuoto)$;
\item we denote by $X_{123}$ the blow-up of $X_{0}$ along the closed
$G$-orbit;
\item we denote by $X_{123,ij}$ the blow-up of $X_{123}$ along
$(cone(e_{i} , e_{j}),\vuoto)$;
\item we denote by $X_{123,i}$ the blow-up of $X_{123}$ along
$(cone(e_{i},e_{1}+e_{2}+e_{3}),\vuoto)$.
\end{itemize}
By \cite{Ru09} Corollary 2.1, if $(G,\grt)= (G_1,\grt)\times (G_2,\grt)$ and $X_0$ is wonderful, then $H=H_1\times H_2$, where  $H_i:=H\cap G_i$. Given an 1-rank symmetric space $G_i/H_i$, let $\psi_{i}(r):=-2\rho+2\rho_0-re_i^*$, $m_i:=max\{r: \psi_i(r)(\gra_i^\cech)<0\}$ and $\bar{m}_i:=max\{r: \psi_i(r)(\gra_i^\cech)\leq0\}$. In Figure \ref{fig rank1
weight} are written the value of $m_i$ and $\bar{m}_i$ for the various $G_i/H_i$. Moreover, we indicate by $e$ (resp. by $h$) when $G_i/H_i$ is exceptional (resp. hermitian non-exceptional).  By an explicitly analysis we can prove the following theorems.

\begin{figure} \label{fig rank1 weight}\begin{changemargin}{1.6cm}{1cm}
\begin{tabular}{|l |l|l|l|l|}\hline

$G/H$&$\theta$& $\bar{m}_1$&$m_1$&\\
\hline \hline

$PSL_{2}$&$A_{1}$& $1$&$0$&\\\hline
$SL_{2}$&$A_{1}$& $2$&$1$&\\\hline
$SL_{2}/N(SO_{2})$&AI& $ 0$& $0$&h\\\hline
$SL_{2}/SO_{2}$&AI& $ 1$& $0$&h\\\hline
$SL_{2}/N(Sp_{2})$&AII&  $2$& $1$&\\\hline
$SL_{2}/Sp_{2}$&AII&  $4$& $3$&\\\hline
$SL_{n}/S(GL_{1}\times GL_{n-1})$&AIV& $ n$& $ n-1$&e\\\hline
$SO_{2n+1}/N(SO_{1}\times SO_{2n})$&BII&  $ n-1$&$ n-1$&\\\hline
$SO_{2n+1}/(SO_{1}\times SO_{2n})$&BII&  $ 2n-1$&$ 2n-2$&\\\hline
$Sp_{2n}/(Sp_{2}\times Sp_{2n-2})$&CII& $2n-1$ & $ 2n-2$&\\\hline
$SO_{2n}/S( O_{1}\times  O_{2n-1})$&DII&$n-1$&$ n-2$&\\\hline
$SO_{2n}/ (SO_{1}\times SO_{2n-1})$&DII&$2n-2$&$2n-3$&\\\hline
$F_{4}/B_{4}$&FII&$11$ & $10$&\\\hline

\end{tabular}\end{changemargin}
\caption{Weights of rank  1 symmetric spaces}\label{fig rank1
weight}
\end{figure}

\begin{figure}\label{fig rank3 indecom}
\begin{changemargin}{0.7cm}{1cm}
\begin{tabular}{|l|c|l|}

\hline
$G/H$& $\theta$&   $X$\\
\hline \hline

$PSL_{4}$&$A_{3}$& $X_{0}$, $X_{13}$\\

\hline $SO_{7}$&$B_{3}$& $X_{0}$\\ \hline

$PSp_{6}$&$C_{3}$& $X_{0}$, $X_{13}$\\ \hline

$Sp_{6}$&$C_{3}$& $X_{0}$, $X_{13}$\\ \hline

$SL_{4}/N (SO_{4})$&AI& $X_{0}$\\ \hline

$SL_{8}/N (Sp_{8})$&AII&$X_{0}$, $X_{12}$, $X_{13}$, $X_{23}$, $X_{123,13}$\\\hline

 $ SL_{n+1}/S(L_{3}\times L_{n-2})$, $n\geq6$ &AIII&
$X_{0}$, $X_{13}$\\\hline

$SL_{6}/N (S(L_{3}\times L_{3}))$&AIII& $X_{0}$\\ \hline

$SL_{6}/S(L_{3}\times L_{3})$&AIII& $X_{0}$\\ \hline

$SO_{2n+1}/N (SO_{3}\times SO_{2n-2})$&BI&$X_{0}$\\\hline

$Sp_{6}/N (GL_{3})$&CI&$\nexists$\\\hline

$Sp_{6}/GL_{3}$&CI&$X_0$\\\hline

$Sp_{12}/N (Sp_{6}\times Sp_{6})$&CII& $X_{0}$, $X_{12}$,
$X_{13}$, $X_{23}$\\ \hline

$Sp_{12}/(Sp_{6}\times Sp_{6})$&CII& $X_{0}$, $X_{12}$, $X_{13}$, $X_{23}$,
$X_{123,13}$\\ \hline

$Sp_{2n}/(Sp_{6}\times Sp_{2n-6})$, $n>6$ &CII&$X_{0}$, $X_{12}$, $X_{13}$,
$X_{23}$, $X_{123,13}$\\
\hline

$SO_{2n }/N (SO_{3}\times SO_{2n-3})$&DI& $X_{0}$\\ \hline

$SO_{12}/N (GL_{6})$&DIII& $X_{0}$, $X_{12}$\\ \hline

$SO_{12}/N (GL_{6})$&DIII& $X_{0}$, $X_{12}$\\ \hline

$SO_{14}/GL_{7}$&DIII& $X_{0}$, $X_{12}$, $X_{13}$, $X_{23}$,
$X_{123,13}$\\ \hline

$E_{7}/N (E_{6}\times\mathbb{C}^{*})$&EVII& $X_{0}$, $X_{12}$\\ \hline

$E_{7}/ (E_{6}\times \mathbb{C}^{*})$&EVII& $X_{0}$, $X_{12}$\\
\hline

\end{tabular}\end{changemargin}
\caption{Fano indecomposable symmetric varieties of rank
3}\label{fig rank3 indecom}
\end{figure}

\begin{thm} Let $G/H$ be  an indecomposable  symmetric space of
rank 3 such that its standard compactification $X_0$ is wonderful. If $X$ is a
smooth Fano compactification of $G/H$ obtained  from $X_0$
by a sequence of blow-ups along $G$-subvarieties, then it is
$X_{0}$, $X_{12}$, $X_{13}$, $X_{23}$ or $X_{123,13}$. More
precisely, the Fano ones are those which appears in Figure
\ref{fig rank3 indecom}.

\end{thm}

\begin{figure}\label{fig rank3 decom2+1}\begin{changemargin}{0.6cm}{1cm}
\begin{tabular}{|l|c|l|}

\hline
$G_{2}/H_{2}$& $m_1$&   $X$\\
\hline \hline

$ PSL_{3}$&$-$& $X_{0}$, $X_{12}$, $X_{13}$, $X_{23}$\\

\hline $SO_{5}$&$-$& $X_{0}$, $X_{12}$\\

\hline $Spin_{5}$&$-$& $X_{0}$, $X_{12}$, $X_{13}$, $X_{23}$\\

\hline $G_{2}$&$-$& $X_{0}$, $X_{13}$\\

\hline $SL_{3}/N (SO_{3})$&$-$& $X_{0}$\\

\hline$SL_{6}/N (Sp_{6})$&$-$&$X_{0}$, $X_{12}$, $X_{13}$, $X_{23}$, $X_{123,23}$\\ \cline{2-3}
&$m_1\geq1$&$X_{123}$, $X_{123,12}$, $X_{123,13}$ \\

\hline  $SL_{n+1}/S(GL_{2}\times
GL_{n-1})$&$-$&$X_{0}$, $X_{12}$, $X_{13}$,
$X_{23}$\\\cline{2-3}
($n\geq5$)&$m_1\geq1$&$X_{123,13}$\\

\hline  $SL_{5}/S(GL_{2}\times
GL_{3})$&$-$&$X_{0}$, $X_{12}$, $X_{13}$,
$X_{23}$\\

\hline $SL_{4}/N (S(GL_{2}\times GL_{2}))$&$-$& $X_{0}$\\

\hline $SL_{4}/S(GL_{2}\times GL_{2})$&$-$& $X_{0}$, $X_{12}$\\

\hline $SO_{5}/S(O_{2}\times O_{3})$&$-$& $X_{0}$\\

\hline $SO_{5}/  (SO_{2}\times SO_{3})$ &$-$&  $X_{0}$\\

\hline $SO_{2n+1}/S(O_{2}\times O_{2n-1})$ &$-$& $X_{0}$, $X_{13}$\\

\hline $SO_{2n+1}/  (SO_{2}\times SO_{2n-1})$&$-$& $X_{0}$, $X_{13}$\\

\hline   $Sp_{2n}/Sp_{4}\times Sp_{2n-4}$&$-$&
$X_{0}$, $X_{12}$, $X_{13}$, $X_{23}$, $X_{123,23}$\\ \cline{2-3}
$(n\geq5)$&$m_1\geq1$&$X_{123}$, $X_{123,12}$, $X_{123,13}$\\ \cline{2-3}
&$m_1\geq2$&$X_{123,2}$\\

\hline $Sp_{8}/N (Sp_{4}\times Sp_{4})$&$-$& $X_{0}$,
$X_{12}$, $X_{13}$, $X_{23}$\\

\hline $Sp_{8}/Sp_{4}\times Sp_{4}$&$-$& $X_{0}$, $X_{12}$,
$X_{13}$, $X_{23}$, $X_{123,23}$\\ \cline{2-3}
&$m_1\geq1$&$X_{123}$, $X_{123,12}$, $X_{123,13}$\\

\hline $SO_{2n}/S ( O_{2}\times  O_{2n-2})$&$-$& $X_{0}$, $X_{13}$\\

\hline $SO_{2n}/SO_{2}\times SO_{2n-2}$&$-$&$X_{0}$, $X_{13}$\\

\hline $SO_{8}/N (GL_{4})$&$-$& $X_{0}$, $X_{12}$\\

\hline $SO_{8}/GL_{4}$&$-$&$X_{0}$, $X_{12}$\\

\hline $SO_{10}/ GL_{5} $&$-$& $X_{0}$, $X_{12}$, $X_{13}$,
$X_{23}$, $X_{123,23}$\\ \cline{2-3}
&$m_1\geq1$&$X_{123}$, $X_{123,12}$, $X_{123,13}$\\

\hline $E_{6}/D_{5}\times \mathbb{C}^{*}$&$-$& $X_{0}$,
$X_{12}$, $X_{13}$, $X_{23}$, $X_{123,23}$\\ \cline{2-3}
&$m_1\geq1$&$X_{123}$, $X_{123,12}$, $X_{123,13}$, $X_{123, 1}$\\ \cline{2-3}
&$m_1\geq2$& $X_{123,2}$, $X_{123,3}$\\

\hline $E_{6}/N (F_{4})$&$-$& $X_{0}$, $X_{12}$,
$X_{13}$, $X_{23}$, $X_{123,23}$\\ \cline{2-3}
&$m_1\geq1$& $X_{123}$, $X_{123,12}$, $X_{123,13}$, $X_{123, 1}$ \\ \cline{2-3}
&$m_1\geq2$& $X_{123,2}$, $X_{123,3}$\\

\hline $G_{2}/(A_{1}\times A_{1})$&-& $\nexists$ \\\hline

\end{tabular}\end{changemargin}
\caption{Fano decomposable symmetric varieties of rank 3}\label{fig
rank3 decom2+1}
\end{figure}

\

\begin{thm} Let $G/H$ be a symmetric space such that $X_0$ is
wonderful and such that $(G,\theta)=(G_{1},\theta)\times (G_{2},\theta)$ with
$rank\,G_{i}/G_{i}^{\theta}=i$.   If $X$ is a smooth Fano compactification of
$G/H$ obtained from $X_0$  by a sequence of blow-ups
along $G$-subvarieties, then it is $X_{0}$, $X_{12}$, $X_{13}$,
$X_{23}$, $X_{123}$, $X_{123,12}$, $X_{123,13}$, $X_{123,23}$,
$X_{123,1}$, $X_{123,12}$ or $X_{123,3}$. More precisely, the
classification of such varieties is as in Figure \ref{fig rank3
decom2+1}. In the second column, we indicate the conditions on
$m_1 $ so that $X$ is Fano.

\end{thm}

\

\begin{thm} Let $G/H$ be a symmetric space such that $X_0$   is
wonderful. Suppose that  $(G,\theta)=(G_{1},\theta)\times
(G_{2},\theta)\times(G_{3},\theta) $ with
$rank\,G_{i}/G_{i}^{\theta}=1$ and let $x_{r}$ be the number of factors $G_i$
such that $\psi_{i}(r) $ is antidominant and regular. If $X$ is a
smooth Fano compactification of $G/H$ obtained from $X_0$
by a sequence of blow-ups along $G$-subvarieties, then it is
$X_{0}$, $X_{12}$, $X_{13}$, $X_{23}$, $X_{123}$, $X_{123,12}$,
$X_{123,13}$, $X_{123,23}$, $X_{123,1}$, $X_{123,2}$ or
$X_{123,3}$.   More precisely, we have the following classification
(depending on $G/H$):

\begin{itemize}
\item If $x_{1}\leq1$,  then the smooth  Fano compactifications of $G/H$ are
$X_{0}$, $X_{12}$, $X_{13} $ and
$X_{23} $. In particular, there are four of them.
\item If $x_{1}=2$, then there are five Fano varieties. Let $i<j$ be the indices such that $\psi_{i}(1)$ and $\psi_{j}(1)$
are anti-dominant and regular.  The smooth Fano compactifications
of $G/H$ are  $X_{0} $, $X_{12}$,
$X_{13}$, $X_{23}$ and $X_{123,ij}$.
\item If $(x_{1},x_{2})$ is equal to $(3,0)$ or to $(3,1)$, then there are eight Fano varieties:
$X_{0}$, $X_{12}$, $X_{13}$, $X_{23}$, $X_{123}$, $X_{123,12}$,
$X_{123,13}$ and $X_{123,23}$.
\item If $(x_{1},x_{2})=(3,2)$, then there are nine Fano varieties.
Suppose that the $\psi_{j}(2)$ with $j\neq i$ are
anti-dominant and regular, then  the Fano
varieties are: $X_{0}$, $X_{12}$, $X_{13}$, $X_{23}$, $X_{123}$,
$X_{123,12}$, $X_{123,13}$, $X_{123,23}$ and $X_{123,i}$.
\item If $(x_{1},x_{2})=(3,3)$, then there are eleven Fano varieties:
$X_{0}$, $X_{12}$, $X_{13}$, $X_{23}$, $X_{123}$, $X_{123,12}$,
$X_{123,13}$, $X_{123,23}$, $X_{123,1}$, $X_{123,2}$ and
$X_{123,3}$.
\end{itemize}
\end{thm}

Remark that the first part of the proof holds in a more general contest.

\begin{cor}
Given any wonderful $G$-variety $X_0$ of rank 3 (even non symmetric) and any Fano variety $X$ obtained from $X_0$ by a succession of blow-ups along $G$-stable subvarieties, then $X$ is one of the following varieties: $X_0$, $X_{12}$, $X_{13}$, $X_{23}$, $X_{123}$, $X_{123,12}$, $X_{123,13}$, $X_{123,23}$, $X_{123,1}$, $X_{123,2}$ or $X_{123,3}$.
\end{cor}

\section{(Quasi) Fano symmetric varieties of rank 2}

\subsection{Fano symmetric varieties}
In this section we consider the quasi-Fano locally factorial  symmetric varieties with rank 2 (and
$G$ only reductive). Remark that to individuate univocally a projective 2-rank symmetric variety with $\rho$ injective over $\varrho^{-1}(\varrho(D(X))$ is sufficient to give $\grD[p]$.  We begin classifying the Fano varieties with  $G$ semi-simple.  First, we   consider    two special cases: i)  $G/H$ is indecomposable, while $X$ is neither simple nor toroidal; ii) $X$ is toroidal.

\begin{lem}\label{lem fano indecomposable non toroidal} Let $X$ be a locally factorial projective symmetric variety. Suppose that $G/H$ is indecomposable and that $X$ is neither simple nor toroidal, then
$R_{G,\grt}=A_2$, $H=G^\grt$ and $\calF(X)$ contains $(cone( -\gra_1^\cech-\gra_2^\cech),\vuoto)$.
\end{lem}

{\em Proof.}   We do a case-to-case analysis.
1) Suppose   $R_{G,\grt}=A_2$ and $H=N(G^\grt)$. Then there is   $cone(\gra_1^\cech,v)$ in $\grD$, with $v:=-\frac{x}{3}\gra^\cech_1-\frac{y}{3}\gra^\cech_2\in int(C^-)$ primitive. We can write     $-\gro^\cech_2
$ as a positive integral combination of $\gra^\cech_1$ and $v$; thus   $y\in\{1,2\}$. Moreover, $0<(-v,\gra_2)= -\frac{x}{3}+\frac{2y}{3}$, so $0<x<2y$.
If $y=1$, then 
$v=-\frac{1}{3}(\gro^\cech_1+\gro^\cech_2)\notin\chi_*(S)$. If $y=2$, we have three possibilities: 
i) $v= -\gro^\cech_2 $ which is not regular; 
ii) $v=-\frac{2}{3}(\gro^\cech_1+\gro^\cech_2)\notin \chi_*(S)$;
iii) $v=- \gra^\cech_1-\frac{2}{3}\gra^\cech_2\notin \chi_*(S)$.

2) Suppose   $R_{G,\grt}=A_2$ and $H=G^\grt$. Then there is   $cone(\gra_1^\cech,v)$ in  $\grD$, with $v:=-x\gra^\cech_1-y\gra^\cech_2$ primitive. Hence $y=1$ because $\{\gra_1^\cech,v\}$ is a basis of $\chi_*(S)$.  Moreover, $0<(-v,\gra_2)= -x+2$, so 
$v=-\gra^\cech_1-\gra^\cech_2$ as in the statement.

3) Suppose  that $R_{G,\grt}=B_2$ and that $cone(\gra_1^\cech,v)\in \grD$, with $v:=-x\gra^\cech_1-\frac{y}{2}\gra^\cech_2$ primitive.  We can write $-\gro^\cech_2=-\gra^\cech_1-\gra^\cech_2$ as a positive integral combination of $\gra^\cech_1$ and $v$; thus $y\in\{1,2\}$. Moreover, $0<(-v,\gra_2)=-x+y$, so $0<x<y$. Therefore
$v= -\gro^\cech_2 $ is not regular.

4) Suppose  that $R_{G,\grt}=B_2$ and  that $cone(\gra_2^\cech,v)\in \grD$. If $H=N(G^\grt)$, then  $\{\gra_2^\cech,v\}$ cannot be a basis of $\chi_*(S)=\mZ\gra_1^\cech\oplus\mZ\frac{\gra_2^\cech}{2}$.

5) Suppose   $R_{G,\grt}=B_2$ and $H= G^\grt $. Suppose also that $cone(\gra_2^\cech,v)\in\grD$, with $v:=-x\gra^\cech_1-y\gra^\cech_2$ primitive. Then $x=1$ because $\{\gra_2^\cech,v\}$ is a basis of $\chi_*(S)$.   Moreover, $0<(-v,\gra_1)= 2-2y$, so    $0<y<1$.

6) Suppose   $R_{G,\grt}=BC_2$ and     $cone(\gra_1^\cech,v)\in \grD$, with $v:=-x\gra^\cech_1-y\gra^\cech_2$ primitive. As before $y=1$.   Moreover, $0<(-v,\gra_2)= -x+1$, so $0<x<1$.

7) Suppose    $R_{G,\grt}=BC_2$ and    $cone(\gra_2^\cech,v)\in \grD$, with $v:=-x\gra^\cech_1-y\gra^\cech_2$ primitive. As before $x=1$ and   $0<(-v,\gra_1)= 2-y$, so 
$v=-\gro^\cech_1$ is not regular.

8) Suppose    $R_{G,\grt}=G_2$ and      $cone(\gra_1^\cech,v)\in \grD$, with $v:=-x\gra^\cech_1-y\gra^\cech_2$ primitive. As before $y=1$ and   $0<(-v,\gra_2)\leq -3x+2$, so $0<x<\frac{2}{3}$.

9) Suppose   $R_{G,\grt}=G_2$ and     $cone(\gra_2^\cech,v)\in \grD$, with $v:=-x\gra^\cech_1-y\gra^\cech_2$ primitive. As before $x=1$.   Moreover, $0<(-v,\gra_1)= 2-y$ and $0<(-v,\gra_2)= -3+2y$; thus $3<2y<4$. $\square$

In the next two lemmas, we do not make any hypothesis on the  regularity of $X$.

\begin{lem}\label{numero 1 coni X fano rango 2} Let $X$ be a Fano symmetric variety with $G$ semisimple.
Then $\calF(X)$ contains  at most three  colored 1-cones.  \end{lem}

{\em Proof.} Suppose by contradiction that $\mathcal{F}(X)$ contains $(\sigma_i:=cone(v_i,v_{i+1}),\vuoto)$, with  $i=1,2,3$ and    $v_1,v_2,v_3,v_4$ primitive. We can write $v_2$ as a positive linear combination $xv_{1}+yv_{4}$ of $v_1$ and $v_4$. Then $k_{\sigma_2}(v_{1})<k(v_{1})=1$, $k_{\sigma_2}(v_{4})<k(v_{4})=1$ and $1=k_{\grs_2}(v_2)=xk_{\grs_2}(v_1)+yk_{\grs_2}(v_4)\leq0$. $\square$

\begin{lem}\label{lem fano toroidal} Let $X$ be a Fano non-simple  toroidal  symmetric variety. Then $X_0$ is smooth and $\grD[p]$ is either $\{ e_1,e_1+re_2,e_2\}$ or  $\{ e_1,  re_1+e_2,e_2\}$. These varieties are smooth if and only if $r=1$.
\end{lem}

{\em Proof.} By the previous lemma, we have $\grD[p]=\{e_1,e_2,v\}$ for an appropriate $v$.  First
suppose that   $X_0$ is smooth, i.e.
$\chi_{*}(S)=\mathbb{Z}e_{1}\oplus \mathbb{Z} e_{2}$, and write
$v=x_{1}e_{1}+x_{2}e_{2}$. For each $i$, let $\sigma_{i}=cone(e_{i},v)$ and $\{i,i^c\}=\{1,2\}$, so   $k_{\sigma_{i}}$ is
$e_{i}^{*}+\frac{1-x_{i}}{x_{i^{c}}}e_{i^{c}}^{*}$.    If  $x_{1}\geq2 $ and $x_{2}\geq2$, then
$\frac{1-x_{1}}{x_{2}}$ and $ \frac{1-x_{2}}{x_{1}}$ are strictly negative integers, so    $x_{1}\geq x_{2}+1\geq x_{1}+2$, a contradiction.

Suppose now that $X_0$ is singular. Then $R_{G,\theta}$ is
either $A_{2}$ or $A_{1}\times A_{1}$. In the first case,
$e_{i}=-3\omega_{i}^{\vee}$ and the strictly positive integer
$k(-\omega_{1}^{\vee}-\omega_{2}^{\vee}) $
is lesser than
$\frac{1}{3}k(-3\omega^{\vee}_{1})+\frac{1}{3}k(-3\omega^{\vee}_{2})= \frac{2}{3}$,
a contradiction.

Finally suppose $R_{G,\theta}=A_{1}\times A_{1}$ and
$\chi_{*}(S)=\mathbb{Z}2\omega_{2}^{\vee}\oplus
\mathbb{Z}(\omega_{1}^{\vee}+\omega_{2}^{\vee})$. Then, there is $i$ such that  $-\omega_{1}^{\vee}-\omega_{2}^{\vee} =av+b(-2\gro_{i}^\cech)$ with $a,b\geq0$. The
integer $a+b=k(-\omega_{1}-\omega_{2})$ is strictly lesser than
$\frac{1}{2}k(-2\omega_{1})+\frac{1}{2}k(-2\omega_{2}) =1$, a contradiction.  
$\square$

Remark that the previous two lemmas apply also to a toroidal symmetric variety with $-K_X$ ample, $|\calF(X)|$ convex and generated by a basis of $\chi_*(S)$ (without supposing $X$ compact). Now, we   state the main result of this section.

\begin{thm} Let $G/H$ be a  symmetric space of rank 2 (with $G$ semisimple). \label{fano rank2}

\begin{itemize}
\item If a (projective) symmetric variety $X$ is Fano then $\rho^{-1}(\rho(D(X)))=D(X)$. If moreover $X$ is locally factorial, then $\sharp\rho(D(X))=\sharp\rho^{-1}(\rho(D(X)))$.

\item   If $R_{G,\theta}$ is irreducible and $X$ is a Fano locally factorial equivariant compactification of $G/H$, we have exactly the following
possibilities for $\grD[p]$:
\begin{enumerate}
\item $\Delta[p]=\{e_{1},\alpha_{1}^{\vee}\}$  if $\sharp\rho^{-1}(\gra^\cech_1)=1$ and  $R_{G,\theta}$ is not $G_{2}$;
\item $\Delta[p]=\{e_{2},\alpha_{2}^{\vee}\}$ if $\sharp\rho^{-1}(\gra^\cech_2)=1$ and $R_{G,\theta}$ is not $B_{2}$;
\item $\Delta[p]=\{e_{2},\alpha_{2}^{\vee}\}$ if   $R_{G,\theta}=B_{2}$ and $H=G^\grt$;
\item    $\Delta[p]=\{-\omega_{1}^{\vee},-\omega_{2}^{\vee}\}$ if $G/H\neq G_2/(SL_2\times SL_2)$ and $H=N(G^{\theta})$;
\item $\Delta[p]=\{-2\omega_{1}^{\vee},-\omega_{2}^{\vee}\}$ if $R_{G,\theta}=B_{2}$ and $H=G^{\theta}$;
\item $\Delta[p]=\{\alpha_{1}^{\vee},\alpha_{2}^{\vee},-\alpha_{1}^{\vee}-\alpha_{2}^{\vee}\}$ if  $H=G^{\theta}$ and $R_{G,\theta}=A_{2}$;
\item $\grD[p]=\{ \alpha_{1}^{\vee},-\gro_{1}^{\vee}-\gro_{2}^{\vee},-3\gro_{1}^{\vee}\}$ (and  $\grD[p]=\{\alpha_{2}^{\vee},-\gro_{1}^{\vee}-\gro_{2}^{\vee},-3\gro_{2}^{\vee}\}$) if $H=G^{\theta}$ and  $R_{G,\theta}=A_{2}$;

\item $\grD[p]=\{\alpha_{1}^{\vee},-\gro_{1}^{\vee}-\gro_{2}^{\vee},-4\gro_{1}^{\vee}-\gro_{2}^{\vee}, -3\gro_{1}^{\vee} \}$
(and   $\grD[p]=\{\alpha_{2}^{\vee},-\gro_{1}^{\vee}-\gro_{2}^{\vee},- \gro_{1}^{\vee}-4\gro_{2}^{\vee}, -3\gro_{2}^{\vee} \}$) if $H=G^{\theta}$,
 $R_{G,\theta}=A_{2}$ and $\grt\neq-id$ over $\chi_*(T)$;

\item $\Delta[p]=\{e_{1},e_{2}, e_{1}+e_{2}\}$ if
$\chi_{*}(S)=\mathbb{Z}e_{1}\oplus\mathbb{Z} e_{2} $ and $-2\rho +2\rho_0+e^*_{i}\in int(C^-)$ for each $i$.
\end{enumerate}
The previous varieties are singular in the following cases:
\begin{itemize}
\item $\Delta[p]=\{e_{i},\alpha_{i}^{\vee}\}$ and $i=1,2$  if   $R_{G,\theta}=A_2$;
\item $\Delta[p]=\{e_{1},\alpha_{1}^{\vee}\}$  if $R_{G,\theta}=B_2$, $H=G^\grt$  (and $\sharp D(G/H)=2$);
\item $\Delta[p]=\{e_{2},\alpha_{2}^{\vee}\}$  if $R_{G,\theta}=BC_2$ (and $\sharp D(G/H)=2$).
\end{itemize}

\item Suppose $(G,\theta)=(G_{1},\theta)\times(G_{2},\theta)$ and let   $m_i$ be as in section \S\ref{sez rk3}.   If $\chi_*(S)$ $=\mZ e_1\oplus \mZ e_2$, we have the following locally factorial Fano varieties: 
\begin{enumerate}
\item $\grD[p]=\{e_1,e_2 \}$;
\item $\grD[p]=\{e_1,e_2,e_{1}+e_{2}\}$;
\item $\grD[p]=\{\alpha_{1}^{\vee},-\alpha_{1}^{\vee},-r\alpha_{1}^{\vee}+e _{2}\}$ if $ r\leq m_2+1$ and $\sharp\rho^{-1}(\gra^\cech_1)=1$;
\item $\grD[p]=\{\alpha_{2}^{\vee},-\alpha_{2}^{\vee},e_{1} -r\alpha_{2}^{\vee}\}$ if $ r\leq m_1+1$ and $\sharp\rho^{-1}(\gra^\cech_2)=1$;

\item $\grD[p]=\{\alpha_{1}^{\vee},-\alpha_{1}^{\vee},-r\alpha_{1}^{\vee}+e _{2},-(r+1)\alpha_{1}^{\vee}+e _{2}\}$ if \ $ r\leq m_2 $ \ and $\sharp\rho^{-1}(\gra^\cech_1)=1$;
\item $\grD[p]=\{\alpha_{2}^{\vee},-\alpha_{2}^{\vee}, e_{1} -r\alpha_{2}^{\vee}, e_{1} -(r+1)\alpha_{2}^{\vee}\}$ \ if \ $ r\leq m_1 $ \ and $\sharp\rho^{-1}(\gra^\cech_2)=1$;

\item $\grD[p]=\{\alpha_{1}^{\vee},\alpha_{2}^{\vee},-\alpha_{1}^{\vee}-\alpha_{2}^{\vee}\}$ if $\sharp D(G/H)=2$.

\end{enumerate}
Only the first two are smooth.

\item If $G/H$ is decomposable but   $\chi_*(S)=
    \mathbb{Z}(\frac{1}{2}\alpha_{1}^{\vee}+\frac{1}{2}\alpha_{2}^{\vee})
    \oplus\mathbb{Z}\alpha_{2}^{\vee}$,  we have the following locally factorial  Fano compactifications of $G/H$:
\begin{enumerate}

\item $\grD[p]= \{\alpha_{1}^{\vee},\alpha_{2}^{\vee}, -\frac{1}{2}\alpha_{1}^{\vee}-\frac{1}{2}\alpha_{2}^{\vee}\}$;
\item $\grD[p]=
    \{\alpha_{1}^{\vee},-\alpha^{\vee}_{1}, -\frac{2r+1}{2}\alpha_{1}^{\vee}-\frac{1}{2}\alpha_{2}^{\vee}\}$ if $0\leq r\leq \frac{m_2+1}{2}$;
\item $\grD[p]= \{\alpha_{2}^{\vee},-\alpha^{\vee}_{2}, -\frac{1}{2}\alpha_{1}^{\vee}-\frac{2r+1}{2}\alpha_{2}^{\vee}\}$ if $0\leq r\leq \frac{m_1+1}{2}$;
\item $\grD[p]=
    \{\alpha_{1}^{\vee},-\alpha^{\vee}_{1}, -\frac{2r+1}{2}\alpha_{1}^{\vee}-\frac{1}{2}\alpha_{2}^{\vee},-\frac{2r+3}{2}\alpha_{1}^{\vee}-\frac{1}{2}\alpha_{2}^{\vee}\}$ if $0\leq r\leq \frac{m_2-1}{2}$;
\item $\grD[p]= \{\alpha_{2}^{\vee},-\alpha^{\vee}_{2}, -\frac{1}{2}\alpha_{1}^{\vee}-\frac{2r+1}{2}\alpha_{2}^{\vee},-\frac{1}{2}\alpha_{1}^{\vee}-\frac{2r+3}{2}\alpha_{2}^{\vee}\}$ if $0\leq r\leq \frac{ m_1-1}{2}$.

\end{enumerate}
The first is smooth, while the other ones are smooth if and only if $r=0$.

\end{itemize}
\end{thm}

Remark that, supposing $R_{G,\grt}=A_1\times A_1 $,   the $m_i$ have the same value if $\chi_*(S)=    \mathbb{Z}(\gro_{1}^{\vee}+\gro_{2}^{\vee}) \oplus\mathbb{Z}2\gro_{2}^{\vee}$ or if   $H= G^\grt $. Moreover, $\sharp D(G/H)=2$ if $\chi_*(S)=    \mathbb{Z}(\gro_{1}^{\vee}+\gro_{2}^{\vee}) \oplus\mathbb{Z}2\gro_{2}^{\vee}$, because any $n\in H  \setminus G^\grt$ exchanges  the   colors of $G/G^\grt$ associated to the   same coroot.

{\em Proof.} The first property holds because $-K_X$ is linearly equivalent to a $G$-stable divisor. Let $a_i:=(2\rho-2\rho_0)(\gra_i^\cech)$ for each $i$; recall that $a_1,a_2\geq1$ (see \S\ref{sez standard}). In the following we always suppose $\sharp \rho(D(X))=\sharp\rho^{-1}(\rho(D(X)))$.

I) Suppose $X$  simple.  If $X$ is toroidal (i.e. $X=X_0$) have to exclude the following cases by \S\ref{sez standard}: i) $G/H=G_2/(SL_2\times SL_2)$;  ii) $R_{G,\grt}=A_2$ and $H=G^\grt$; iii) $R_{G,\grt}=A_1\times A_1$ and
$\chi_*(S)=\mZ2\gro_1^\cech\oplus\mZ(\gro_1^\cech+\gro_2^\cech)$. In the first case $X_0$ is not Fano, while in the last two cases  $X_0$ is not smooth.  If $\rho(D(X))=\{\gra^\cech_1\}$, then $\grt$ is indecomposable because $C(X)$ is strictly convex. Moreover,  $\grD[p]=\{\gra^\cech_1,e_1=-x_1\gra^\cech_1-x_2\gra^\cech_2\}$ and $k_X=a_1\gro_1-\frac{x_1a_1+1}{x_2}\gro_2$; so $k_X(\gra^\cech_2)\leq0<a_2$.  Thus we have only to verify if $X$ is locally factorial.  We have to exclude two cases: 1) $\grD(2)=\{cone(\gra_1^\cech,-\gro_1^\cech
)\}$ and $R_{G,\grt}=G_2$; 2) $\grD(2)=\{cone(\gra_2^\cech,-\gro_2^\cech
)\}$, $R_{G,\grt}=B_2$ and $H=N(G^\grt$).

II) Suppose now that $X$ is not simple. In the toroidal case $\grD[p]=\{e_1,e_2,e_1+e_2\}$ and $X_0$ is smooth by Lemma \ref{lem fano toroidal}. This variety is Fano if and only if $-2\rho+2\rho_0+e_i^*$ is antidominant and regular for each $i$. III) Assume moreover that $X$ is not toroidal.  Suppose first $\grt$ indecomposable.
Then, by Lemma \ref{lem fano indecomposable non toroidal}, $R_{G,\grt}=A_2$, $H=G^\grt$ and   $\grD(2)$ contains $cone(\alpha_{1}^{\vee},-\alpha_{1}^{\vee}-\alpha_{2}^{\vee})$, up to reindexing. If $\rho(D(X))$ contains also $\alpha_{2}^{\vee}$ then $\grD(2)=\{cone(\alpha_{1}^{\vee},-\alpha_{1}^{\vee}-\alpha_{2}^{\vee}), cone(\alpha_{2}^{\vee},-\alpha_{1}^{\vee}-\alpha_{2}^{\vee})\}$.
Moreover, $\widetilde{k}_X(\gra_1^\cech)=\widetilde{k}_X(\gra^\cech_2)=0$  and $\widetilde{k}_X(-\gra^\cech_1-\gra^\cech_2)\geq1$, so $\widetilde{k}_X$ is strictly convex.

If $\rho(D(X))=\{\gra_1^\cech\}$, then $\grD[p]$ contains $-3\gro_1^\cech
$. Remark that $\{-\alpha_{1}^{\vee}-\alpha_{2}^{\vee},-3\gro_1^\cech\}$ is a basis of $\chi_*(S)$. By the Lemma \ref{lem fano toroidal}   we have two possibilities for $\grD[p]$: $\{ \alpha_{1}^{\vee}, -\alpha_{1}^{\vee}-\alpha_{2}^{\vee},-3\gro_{1}^{\vee}\}$ or $\{ \alpha_{1}^{\vee}, -\gra_{1}^{\vee}-\alpha_{2}^{\vee}, -4\gro_{1}^{\vee}-\gro_{2}^{\vee},-3\gro_{1}^{\vee} \}$.
In the first case there are not other conditions because   the weights of $k$ are  $a_1\gro_1-(a_1+1)\gro_2$ and $-\gro_2$. In the last case we have to impose that $a_1,a_2>1$. Indeed the weights of $k$ are $a_1\gro_1-(a_1+1)\gro_2$, $-\gra_2$ and $-\gro_1+\gro_2$. Moreover $(-\gra_2)(\gra^\cech_1)< a_1$ and $(-\gro_1+\gro_2)(\gra^\cech_2)< a_2$. In $\S\ref{sez standard}$, we have seen that  $a_1\leq1$, $a_2\leq1$  if and only if $\grt=-id$.

IV) Suppose now $\grt$ decomposable.  First suppose $D(X)=D(G/H)$.
Then $\grD(2)$ contains $cone(\gra^\cech_1,x(-m\gra^\cech_1-\gra^\cech_2))$ and
$cone(\gra^\cech_2,x(- \gra^\cech_1-r\gra^\cech_2))$ with $x\in\{1,\frac{1}{2}\}$. Remark that $\mR^+(-m\gra^\cech_1-\gra^\cech_2)=\mR^+(- \gra^\cech_1-\frac{1}{m}\gra^\cech_2)$ and
$x(- \gra^\cech_1-r\gra^\cech_2)\notin cone(\gra^\cech_1,x(-m\gra^\cech_1-\gra^\cech_2))$ so $r\leq\frac{1}{m}$. Therefore  $\grD(2)=\{cone(\gra^\cech_1,x(-\gra^\cech_1-\gra^\cech_2)), cone(\gra^\cech_2,x(- \gra^\cech_1-\gra^\cech_2))\}$.

V) Suppose now that $\rho(D(X))$ contains exactly one coroot, say
$\alpha_{1}^{\vee}$. Suppose $X_0$ smooth and let $\grs_0:=cone(\alpha_{1}^{\vee},-m\gra_1^\cech+e_2)$ be in $\grD(2)$. The case with $X_0$ singular is very similar.
We can apply the Lemma \ref{lem fano toroidal} to the maximal open toroidal subvariety $X'$ of $X$ (whose colored fan has support $cone(-\alpha_{1}^{\vee},-m\gra_1^\cech+e_2)$). There are two possibilities for $\grD_{X'}$: its maximal  cones are either  $\{\grs_1:=cone(-\alpha_{1}^{\vee},-m\gra_1^\cech+e_2)\}$ or $\{\grs_1:= cone(-m\gra_1^\cech+e_2,-(m+1)\gra_1^\cech+e_2),\grs_2:=cone(-\alpha_{1}^{\vee},-(m+1)\gra_1^\cech+e_2)\}$.
In the first case the   unique non-trivial condition is  $k_{\grs_1}( \gra^\cech_2)<a_2$, which  is equivalent to the following one:
$\widetilde{k}_{\grs_1}( \gra^\cech_2)=(\psi_2(m-1))( \gra^\cech_2)<0$. In the second case the   unique non-trivial condition is   $\widetilde{k}_{\grs_2}( \gra^\cech_2)=(\psi_2(m))( \gra^\cech_2)<0$.


Now, we verify the smoothness of such varieties. First,  we explain the conditions for a projective locally factorial symmetric variety $X$ with rank two to be smooth (see \cite{Ru1}, Theorems 2.1 and 2.2). Let $Y$ be a open simple $G$-subvariety of
$X$ whose closed orbit is compact; then the associated colored cone $(C,F)$ has dimension two. Write $C=cone(v_1,v_2)$ and let $C^{\cech}=cone(v_1^*,v_2^*)$ be the dual cone; we take  $v_1$ and $v_2$ primitive. If $X$ is smooth then $\rho$ is injective over $F$ and $\rho(F)$ does not contain any exceptional root. Suppose such conditions verified and let $R'$ the root subsystem of $R_{G,\grt}$ generated by the simple roots $\gra$ such that $\gra^\cech\in\rho(F)$ and $\sharp\rho^{-1}(\gra^\cech)=1$. If there   is not such a root, then $Y$ is smooth. Otherwise $Y$ it is smooth if and only if: i) $R'$ has type $A_1$; ii) up to reindexing, $\frac{1}{2}(2v_1^*-v_2^*)$ is the fundamental weight of $R'$.

Suppose now  that $C= cone(\gra^\cech_1,-r\gra^\cech_1-\gra^\cech_2)$ with  $-r\gra^\cech_1-\gra^\cech_2$ primitive. If   $X$ is Fano, then $\sharp\rho^{-1}(\gra^\cech_1)=1$.
If moreover $\grt$ is indecomposable, then $R_{G,\grt}=A_2$,  $H=G^\grt$ and $r=1$.
Furthermore $R'$ is $A_1$ and  $C^\cech=cone(v^*_1=\gro_1-\gro_2,v^*_2=-\gro_2)$. Hence $\frac{1}{2}(2v^*_1-v^*_2)
=\frac{1}{2}\gra_1$; thus $Y$ is smooth. If $G/G^\grt=G_1/G_1^\grt\times G_2/G_2^\grt$, then  $R'=R_{G_1,\grt}$, $H=G_1^\grt\times (H\cap G_2)$  and $C^\cech=cone(v^*_1=\gro_1-r\gro_2,v^*_2=-\gro_2)$. Hence $R_{G_1,\grt}$ has to be $A_1$ and   $\frac{1}{2}(2v^*_1-v^*_2)=\gro_1+ (\frac{1}{2}-r)\gro_2$ has to be $\frac{1}{2}\gra_1$. Thus $Y$ is not smooth. Suppose now that $R_{G,\grt}=A_1\times A_1$,
$\chi_*(S)=\mZ2\gro_1^\cech\oplus\mZ(\gro^\cech_1+\gro^\cech_2)$  and $C=cone(\gra^\cech_1,-\frac{2r+1}{2}\gra^\cech_1-\frac{1}{2}\gra^\cech_2)$.  Then $R'$ is $A_1$ and $C^\vee=cone(v^*_1=\gro_1-(2r+1)\gro_2,v^*_2=-2\gro_2)$; moreover, $\frac{1}{2}(2v^*_1-v^*_2)=\gro_1-2r\gro_2$; thus $Y$ is smooth if and only if $r=0$. The other cases are similar.
$\square$

Using the Lemma \ref{lem fano toroidal}, one can easily prove the following proposition:

\begin{prop}
Let $G/H$ be a symmetric space of rank 2 with $X_0$ smooth.  If $\grt$ is indecomposable, the Fano toroidal compactifications of $G/H$ are as in figure \ref{fig rank2  toroidal fano indecomposable}. If $\grt$ is decomposable, let $m_i$ be as in \S\ref{sez rk3}. Then,
the Fano toroidal compactifications of $G/H$ are the following ones:
\begin{itemize}
\item $\grD[p]=\{e_1,e_2\}$;
\item $\grD[p]=\{e_1,e_2,re_1+e_2\}$ with $r\leq m_2+1$;
\item $\grD[p]=\{e_1,e_2, e_1+re_2\}$ with $r\leq m_1+1$.
\end{itemize}
\end{prop}

\subsection{Smooth quasi Fano varieties}

Now, we want to classify the smooth (resp.  locally factorial) quasi Fano symmetric varieties
with rank two and $G$ semisimple. A Gorenstein (projective) variety is called
quasi Fano if its anticanonical divisor is big and nef.

\begin{figure}
\begin{changemargin}{-2.5cm}{1cm}
\onehalfspacing
\begin{tabular}{ |l|l|l|l|}

\hline
 $G/H$ &$\theta$& $\Delta[p]$ & $\Delta[p]$ \\ \hline \hline

$PSL_{3}$&$A_{2}$&$\{-\gro^\cech_{1},-\gro^\cech_{2}\}$ & $\{-\gro^\cech_{1},-\gro^\cech_{2}, -\gro^\cech_{1}-\gro^\cech_{2}\}$     \\\hline


$SO_{3}$&$B_{2}$& $\{-\gro^\cech_{1},-\gro^\cech_{2}\}$  &   \\\hline

$Spin_{3}$&$B_{2}$&$\{-2\gro^\cech_{1},-\gro^\cech_{2}\}$& $\{-2\gro^\cech_{1},-\gro^\cech_{2},-2\gro^\cech_{1}-\gro^\cech_{2}\}$
\\\hline

$ G_{2}$&$G_2$& $\{-\gro^\cech_{1},-\gro^\cech_{2}\}$ &   \\ \hline

$SL_{3}/N (SO_{3})$ &$AI$&$\{-\gro^\cech_{1},-\gro^\cech_{2}\}$  &  \\\hline


$SL_{6}/N (Sp_{6})$ &$AII$&$\{-\gro^\cech_{1},-\gro^\cech_{2}\}$ &$\{-\gro^\cech_{1},-\gro^\cech_{2},-2\gro^\cech_{1}-\gro^\cech_{2}\}$ \\\cline{3-4}&&$\{-\gro^\cech_{1},-\gro^\cech_{2}, -\gro^\cech_{1}-\gro^\cech_{2}\}$
 &$\{-\gro^\cech_{1},-\gro^\cech_{2}, -\gro^\cech_{1}-2\gro^\cech_{2}\}$  \\\hline


$ SL_{n+1}/S(GL_{2}\times GL_{n-1})$ &$AIII$&$\{-\gro^\cech_{1},-\gro^\cech_{2}\}$ &$\{-\gro^\cech_{1},-\gro^\cech_{2},-r\gro^\cech_{1} -\gro^\cech_{2}\} $, $ r\leq n-3$     \\\hline

$ SL_{4}/N (S(GL_{2}\times GL_{2})) $&$AIII$ &$\{-\gro^\cech_{1},-\gro^\cech_{2}\}$  &\\\hline

$ SL_{4}/S(GL_{2}\times GL_{2})$& $AIII$& $\{-\gro^\cech_{1},-2\gro^\cech_{2}\}$  &   \\\hline

$ SO_{2n+1}/S (O_{2}\times O_{2n-1}) $&$BI$&$\{-\gro^\cech_{1},-\gro^\cech_{2}\}$  &  \\ \hline

$ SO_{2n+1}/SO_{2}\times SO_{2n-1}$&$BI$&  $\{-2\gro^\cech_{1},-\gro^\cech_{2}\}$ & \\\hline

$ Sp_{2n}/ (Sp_{4}\times Sp_{2n-4})$ &$CII$ &  $\{-\gro^\cech_{1},-\gro^\cech_{2}\}$ &$\{-\gro^\cech_{1},-\gro^\cech_{2},-r\gro^\cech_{1}-\gro^\cech_{2}\}$, $ r\leq 2n-6$
\\\cline{3-4}$n\geq5$&& $\{-\gro^\cech_{1},-\gro^\cech_{2},-\gro^\cech_{1}- 2\gro^\cech_{2}\}$ &   \\\hline

$ Sp_{8}/N(Sp_{4}\times Sp_{4})$ &$CII$& $\{-\gro^\cech_{1},-\gro^\cech_{2}\}$ &$\{-\gro^\cech_{1},-\gro^\cech_{2},-\gro^\cech_{1}-\gro^\cech_{2}\}$ \\
\hline

$Sp_{8}/(Sp_{4}\times Sp_{4})$ &$CII$& $\{-\gro^\cech_{1},-2\gro^\cech_{2}\}$ &$\{-\gro^\cech_{1},-2\gro^\cech_{2}, -2\gro^\cech_{1}-2\gro^\cech_{2}\}$ \\\cline{3-4}&&$\{-\gro^\cech_{1},-2\gro^\cech_{2},- \gro^\cech_{1}-2\gro^\cech_{2}\}$& $\{-\gro^\cech_{1},-2\gro^\cech_{2}, -\gro^\cech_{1}-4\gro^\cech_2\}$  \\
\hline

$SO_{2n}/S(O_{2}\times O_{2n-2})$&$DI$& $\{-\gro^\cech_{1},-\gro^\cech_{2}\}$&  \\\hline

$SO_{2n}/SO_{2}\times SO_{2n-2}$&$DI$&  $\{-2\gro^\cech_{1},-\gro^\cech_{2}\}$ &\\\hline

$ SO_{8}/N (SO_{4})$& $DIII$&  $\{-\gro^\cech_{1},-\gro^\cech_{2}\}$  &\\\hline

$ SO_{8}/SO_{4}$ &$DIII$&  $\{-\gro^\cech_{1},-2\gro^\cech_{2}\}$ &\\\hline

$SO_{10}/SO_{5}$ &$DIII$& $\{-\gro^\cech_{1},-\gro^\cech_{2}\}$  & $\{-\gro^\cech_{1},-\gro^\cech_{2},- \gro^\cech_{1}-\gro^\cech_{2}\}$ \\\cline{3-4}
&& $\{-\gro^\cech_{1},-\gro^\cech_{2},-\gro^\cech_{1}-2\gro^\cech_{2}\}$&   $\{-\gro^\cech_{1},-\gro^\cech_{2},-2\gro^\cech_{1}-\gro^\cech_{2}\}$  \\
\hline

$E_{6}/D_5\times\mC^*$ &$EIII$&$\{-\gro^\cech_{1},-\gro^\cech_{2},-\gro^\cech_{1}-r\gro^\cech_{2}\}$, $r\leq3$   &$\{-\gro^\cech_{1},-\gro^\cech_{2},-r\gro^\cech_{1}-\gro^\cech_{2}\}$, $r\leq 4$   \\ \cline{3-4}
&&$\{-\gro^\cech_{1},-\gro^\cech_{2}\}$&\\\hline

$E_{6}/N (F_{4})$&$EIV$& $\{-\gro^\cech_{1},-\gro^\cech_{2}\}$ &  \\\cline{3-4}
&&$\{-\gro^\cech_{1},-\gro^\cech_{2},-r\gro^\cech_{1}-\gro^\cech_{2}\}$, $r\leq4$   &$\{-\gro^\cech_{1},-\gro^\cech_{2},-\gro^\cech_{1}-r\gro^\cech_{2}\}$, $r\leq4$   \\\hline


$ G_{2}/(A_{1}\times A_{1})$&$G$& $\nexists$ &\\\hline

\end{tabular}\end{changemargin}
\caption{Fano toroidal indecomposable symmetric varieties with rank 2}
\label{fig rank2  toroidal fano indecomposable}
\end{figure}


\begin{thm}\label{simmetriche r2 qfano} Let $G/H$ be a  symmetric space of rank 2 (with $G$ semisimple).

\begin{itemize}
\item The nefness of the anticanonical bundle of a compactification  of $G/H$ depends only by the
fan associated to the colored fan (and not by the whole colored fan).

\item The fans of the locally factorial quasi Fano (but non-Fano) compactifications of an indecomposable symmetric space of rank 2 (with $G$ semisimple) are
whose in Figure  \ref{fig qFano rank2 indecomposable} 
(we have also to request that $\rho$ is injective over $D(X)$). Such a variety is singular if and only if $\rho(D(X))$ contains an exceptional root.

\item If   $(G,\theta)=(G_{1},\theta)\times(G_{2},\theta)$, let   $m_i$ and  $\bar{m}_i$ as in \S\ref{sez rk3}.    Supposing $\chi_*(S)=\mZ e_1\oplus \mZ e_2$, let $v_j(i):=-i\gra_j^\cech+e_{j^c}$ and $w_j(x,y):=-(xy+1)\gra_j^\cech+ye_{j^c}$.     We have the following locally factorial quasi Fano compactifications of $G/H$, which are not Fano (we always suppose $\rho$ injective over $D(X)$):
\begin{enumerate}
\item $\grD[p]=\{\alpha_{1}^{\vee},-\alpha_{1}^{\vee},v_1(r),v_1(r+1),...,v_1(r+s)\}$ if i) $s=0,1$, ii) $r+s\leq \bar{m}_2+1$, iii) either $r+s> m_2+1$  or  $\sharp\rho^{-1}(\gra^\cech_1)=2$;
\item $\grD[p]=\{\alpha_{2}^{\vee},-\alpha_{2}^{\vee},v_2(r),v_2(r+1),...,v_2(r+s)\}$ if i) $s=0,1$, ii) $r+s\leq \bar{m}_1+1$, iii) either $r+s> m_1+1$  or  $\sharp\rho^{-1}(\gra^\cech_2)=2$;
\item $\grD[p]=\{\alpha_{1}^{\vee},-\alpha_{1}^{\vee},v_1(r),v_1(r+1),...,v_1(r+s)\}$ if i) $s\geq2$ and ii) $r+s\leq \bar{m}_2+1$;
\item $\grD[p]=\{\alpha_{2}^{\vee},-\alpha_{2}^{\vee},v_2(r),v_2(r+1),...,v_2(r+s)\}$ if i) $s\geq2$ and ii) $r+s\leq \bar{m}_1+1$;

\item $\grD[p]=\{\alpha_{1}^{\vee},-\alpha_{1}^{\vee}, -r\gra_1^\cech+ e_{2},w_1(r,1),..., w_1(r,s)\}$ if $r\leq \bar{m}_2$ and $2\leq s\leq \bar{m}_1+1$;
\item $\grD[p]=\{\alpha_{2}^{\vee},-\alpha_{2}^{\vee},-r\gra_2^\cech+ e_{1},w_2(r,1),..., w_2(r,s)\}$ if $r\leq \bar{m}_1$ and $2\leq s\leq \bar{m}_2+1$;
\item $\grD[p]=\{\alpha_{1}^{\vee}, \alpha_{2}^{\vee},-\gra^\cech_1-\gra^\cech_2\}$ if $\sharp D(G/H)>2$;

\item $\grD[p]=\{e_{1},e_{1}+e_{2},e_{1}+2e_{2},...,
e_{1}+(s-1)e_{2},e_{1}+se_{2},e_{2})\}$ if $2\leq s\leq \bar{m}_1+1$;
\item $\grD[p]=\{e_{2},e_{1}+e_{2},2e_{1}+e_{2},...,
(s-1)e_{1}+e_{2},se_{1}+e_{2},e_{1})\}$ if $2\leq s\leq \bar{m}_2+1$.

\end{enumerate}

These varieties are smooth if either  they are toroidal  or if, $\ogni \gra^\cech\in\rho(D(X))$,   $\sharp\rho^{-1}(\gra^\cech)=2$ and $2\gra\nin R_{G,\grt}$.

\item If $G/H$ is decomposable but   $\chi_*(S)=\mathbb{Z}(\frac{1}{2}\alpha_{1}^{\vee}+\frac{1}{2}\alpha_{2}^{\vee})
    \oplus\mathbb{Z}\alpha_{2}^{\vee}$, let $v_j(i):=-\frac{2i+1}{2}\gra_j^\cech-\frac{1}{2}\gra_{j^c}^\cech$ and $w_j(x,y):=-\frac{2xy+y+2}{2}\gra_j^\cech-\frac{y}{2}\gra^\cech_{j^c}$.  We have the following  locally factorial quasi Fano compactifications of $G/H$, which are not Fano:
\begin{enumerate}

\item $\grD[p]=\{\alpha_{1}^{\vee},-\alpha_{1}^{\vee},v_1(r),v_1(r+1),...,v_1(r+s)\}$ if i) $r\geq0$, ii) $s\geq2$ and iii) $ r+s\leq \frac{\bar{m}_2+1}{2}$;
\item $\grD[p]=\{\alpha_{2}^{\vee},-\alpha_{2}^{\vee},v_2(r),v_2(r+1),...,v_2(r+s)\}$ if i) $r\geq0$, ii) $s\geq2$  and iii)  $  r+s\leq \frac{\bar{m}_1+1}{2}$;

\item $\grD[p]=\{\alpha_{1}^{\vee},-\alpha_{1}^{\vee},v_1(r),v_1(r+1),...,v_1(r+s)\}$ if i) $r\geq0$, ii) $s=0,1$ and iii) $\frac{m_2+1}{2}< r+s\leq \frac{\bar{m}_2+1}{2}$;
\item $\grD[p]=\{\alpha_{2}^{\vee},-\alpha_{2}^{\vee},v_2(r),v_2(r+1),...,v_2(r+s)\}$ if i) $r\geq0$, ii) $s=0,1$ and iii) $\frac{m_1+1}{2}<  r+s\leq \frac{\bar{m}_1+1}{2}$;

\item $\grD[p]=\{\alpha_{1}^{\vee},-\alpha_{1}^{\vee},-\frac{2r+1}{2}\gra_1^\cech-\frac{1}{2}\gra^\cech_2,w_1(r,1),..., w_1(r,s)\}$ if i) $ r \leq \frac{\bar{m}_2- 1}{2}$ and ii)  $2\leq s\leq \bar{m}_1+1$;
\item $\grD[p]=\{\alpha_{2}^{\vee},-\alpha_{2}^{\vee},-\frac{1}{2}\gra^\cech_1-\frac{2r+1}{2}\gra_2^\cech,w_2(r,1),..., w_2(r,s)\}$ if i) $r \leq \frac{\bar{m}_1- 1}{2}$ and ii)  $2\leq s\leq \bar{m}_2+1$;
\item $\grD[p]=\{-\alpha_{1}^{\vee},-\alpha_{2}^{\vee},-\frac{1}{2}\alpha_{1}^{\vee}-\frac{1}{2}\alpha_{2}^{\vee}\}$.
\end{enumerate}
The last variety is smooth, while the other ones are  smooth if and only if $r=0$.

\end{itemize}
\end{thm}

The idea of the proof is to utilize Theorem \ref{fano rank2} thanks to the following lemma. Remark that, when $\grt$ is indecomposable, if $-K_X$ is nef it is also big (see \S\ref{sez Picard}).

{\em Proof.} The first property holds because the inequalities in the conditions for the nefness of a Cartier divisor  are not strict.  I) First suppose $\rho$ non-injective. We have to consider the varieties whose fan is as in Theorem \ref{fano rank2}, but which are not Fano because    $\sharp\rho^{-1}(\rho(D(X))\neq\sharp\rho(D(X))$. In such cases  $\rho$ has to be injective over $D(X)$, so that $X$ is locally factorial. If $\grt$ is indecomposable we have the following possibilities: i) $R_{G,\grt}=BC_2$ and $\grD(2)=\{cone(\gra^\cech_2,-\gro_2)\}$; ii)  $R_{G,\grt}=B_2$, $H=G^\grt$ and $\grD(2)=\{cone(\gra^\cech_1,-2\gro_1)\}$. If $\grt$ is decomposable, we have to consider all the possibilities listed in Theorem \ref{fano rank2} which correspond to non-toroidal varieties with $\chi_*(S)=\mZ e_1\oplus\mZ e_2$. We have also to substitute the conditions of type $\psi_i(m)<0$ with the corresponding  conditions
$\psi_i(m)\leq0$.

%

\begin{lem}\label{riduzione a fano}
Let $X$ be a projective  symmetric variety   with $-K_{X}$ nef, then there  is a  symmetric
variety $X'$ below $X$ such that the piecewise linear function
associated to $-K_{X'}$ is strictly convex (over the colored fan of $X'$) and coincides with the
function associated to $-K_{X}$. If $X$ is toroidal, we can choose $X'$ toroidal.
\end{lem}

{\em Proof of Lemma \ref{riduzione a fano}}. Let  $\Delta'$ be the fan whose maximal   cones are the maximal cones over which $k_{X}$ is linear.   Given any  cone  $C\in\grD'$, define $F_{C}$ as $\{D\in D(G/H): \rho(D)\in C\}$.   We claim that  $\{(C,F_{C})\}_{C\in \Delta', Int\,C\cap C^-\neq\vuoto}$ is a colored fan  associated to a symmetric variety which satisfies the conditions of the lemma. Remark that $|\grD'|=|\calF(X)|$. We have only  to prove that the maximal cones of $\Delta'$ are strictly convex. Let $C=cone(\alpha_{1}^\cech,...,\alpha_{r}^\cech,-\varpi_{1}^\cech ,...,-\varpi_{s}^\cech)\in \Delta'(l)$ (with $\varpi_{1}^\cech,...,\varpi_{s}^\cech$ dominant) and suppose by contradiction that it contains the line generated by $v=\sum_{i=1}^r a_{i}\alpha_{i}^\cech+\sum_{j=1}^s b_{j}(-\varpi_{j}^\cech)$. Then $C$ contains also $\mR v'$, where  $v'= \sum b_{j}(-\varpi_{j}^\cech)$; indeed $-v'=-v+\sum a_{i}\alpha_{i}^\cech$. Write $-v'=\sum_{\gra\in\ol{R}_{G,\grt}} c_{\gra}\alpha^\cech$; then  $C$ contains all the $ \alpha^\cech$ such that $c_{\gra}\neq0$, because $-v'\in C\cap C^+$ and any spherical weight is a positive rational combination of the simple restricted roots. Thus, if $v'\neq0$, $k_C(-v')=\sum c_{\gra}k_C(\alpha^\cech)\geq0$, while $k_C(v')=\sum b_{j}k_C(-\varpi_{j}^\cech)>0$. Suppose now that $v=\sum_{i=1}^{r} a_{i}\alpha_{i}^\cech$ with $a_{j}\neq0$. Then  $C$ contains  $\mR  \alpha_{j}^\cech$. Write $-\alpha_{j}^\cech=\sum_{i=1}^r a'_{i}\alpha_{i}^\cech+\sum_{i=1}^s b'_{i}(-\varpi_{i}^\cech)$ with positive coefficients, then there is $j_0$ such that  $b_{j_{0}}\neq0$. So $k_C(-\alpha_{j}^\cech)\geq b_{j_{0}}>0$ and $k_C(\alpha_{j}^\cech)>0$, a contradiction. $\square$.

\

II) Suppose  $X$   toroidal and let $X'$ be as in Lemma  \ref{riduzione a fano}. If $X'=X$ and is  simple, then $G/H=G_2/(SL_2\times SL_2)$.
If $X'\neq X$ and is  simple, then it must be singular. Otherwise, given any $w=x_{1}e_{1}+x_{2}e_{2}\in\Delta_X[p]\setminus\Delta_{X'}[p]$,  we have $x_{i}\in\mZ^{>0}$  and
$x_{1}+x_{2}=k(w)=1$. If $R_{G,\grt}=A_2$ and $H=G^\grt$, then $k_X=-\frac{1}{3}(\gro_1+\gro_2)\notin \chi(S)$, so $X$ is not locally factorial. Thus   $R_{G,\grt}=A_{1}\times A_{1}$ and $\chi_{*}(S)= \mathbb{Z}2\gro_{1}^{\vee}\oplus\mathbb{Z}(\gro_{1}^{\vee}+\gro_{2}^{\vee})$.
Let $w=-x_{1}\alpha_{1}^{\vee}-x_{2}\alpha_{2}^{\vee}$ be in $\Delta_X[p]\setminus \Delta_{X'}[p]$; then   $2x_{1},  2x_{2}\in \mZ^{>0}$ and   $x_{1}+x_{2}=k(w)$ is 1. Thus  $\Delta_X[p]=\{-\alpha_{1}^{\vee},-\alpha_{2}^{\vee},-\frac{1}{2}\alpha_{1}^{\vee}-\frac{1}{2}\alpha_{2}^{\vee}\}$. In this case  $-K_X$ is nef but non-ample.

If   $X'$ is not simple, then the standard compactification of $G/H$ has to be smooth because of Lemma \ref{lem fano toroidal}. So it is sufficient to prove the following lemma.

\begin{lem}\label{toroidal qfano rango2}
The smooth toric varieties birationally proper over $\mA^{2}$
with nef anticanonical divisor are, up to isomorphisms,
$\mA^{2}$ and the varieties $Z_{m}$, where $Z_{m}$ is the variety whose fan has maximal
cones $\{cone(e_{1},e_{1}+e_{2}),cone(e_{1}+e_{2},e_{1}+2e_{2}),...,
cone(e_{1}+(m-1)e_{2},e_{1}+me_{2}),cone(e_{1}+me_{2},e_{2})\}$.
\end{lem}

{\em Proof of Lemma \ref{toroidal qfano rango2}.} The piecewise linear $k_{m}$
function associated to the anticanonical bundle of $Z_{m}$
is linear on $cone(e_{1},e_{1}+me_{2})$. 
It is easy to see that thus such function is convex.

Now we  show that, given any $Z$ as in the hypotheses, it is isomorphic to $Z_m$ for an appropriate $m$. Notice, that any smooth toric
variety birationally proper over $\mA^{2}$ is obtained
by a sequence of blow-ups. Thus there is nothing to prove if
$\sharp\Delta[p]<4$. Suppose now that $\sharp\Delta[p]\geq4$; we claim  that  up to isomorphisms, $\Delta$ contains
$cone(e_{1},e_{1}+e_{2})$ and $cone(e_{1}+e_{2},e_{1}+2e_{2})$.
We know that $\tau=\mathbb{R}^{\geq0}(e_{1}+e_{2})$ is contained in
$\Delta$.

First of all, we determine the restrictions of $k$ to the cones
containing $\tau$ and afterwards we will   determine the cones
themselves. Let
$\sigma=cone(e_{1}+e_{2},b_{1}e_{1}+b_{2}e_{2})\in\Delta(2)$ be a
maximal cone containing $\tau$ and write $k_{\sigma}
=x_{1}e^{*}_{1}+x_{2}e^{*}_{2}$, so
$(k_{\sigma})(e_{1}+e_{2})=x_{1}+x_{2}=1$ and
$x_{i}=(k_{\sigma})(e_{i})\leq1$ for each $i$. Hence  $k_{\sigma}$
is  $e^{*}_{i}$ for an appropriate $i$ and $b_{i}= (k_{\sigma})(
b_{1}e_{1}+b_{2}e_{2})=1$.  Because of the non-singularity of $Z$
the only possibilities for $\sigma$ are $cone(e_{1}+e_{2},e_{1})$,
$cone(e_{1}+e_{2},e_{2})$, $cone(e_{1}+e_{2},e_{1}+2e_{2})$ and
$cone(e_{1}+e_{2},2e_{1}+e_{2})$. The fan $\Delta$
does not contain both $cone(e_{1}+e_{2},e_{1}+2e_{2})$ and
$cone(e_{1}+e_{2},2e_{1}+e_{2})$; otherwise
$k_{cone(e_{1}+e_{2},e_{1}+2e_{2})}(2e_{1}+e_{2})=
2>1=k(2e_{1}+e_{2})$. Observe that if $\Delta$
contains $cone(e_{1}+e_{2},2e_{1}+e_{2})$, then $Z$ is isomorphic to
a variety whose fan contains $cone(e_{1}+e_{2},e_{1}+2e_{2})$
by the isomorphism induced by the automorphism of $\chi_*(S)$ that
exchanges $e_{1}$ and $e_{2}$. So the claim is proved.

Because of the non-singularity of $Z$, $\Delta$ contains a cone
$\sigma=cone(e_{1}+me_{2},e_{2})$ for a suitable integer $m$; we
want to show that $Z$ is $Z_{m}$. Let  $Z'$ be the open toric
subvariety of $Z$ whose fan $\Delta'$ is
$\Delta\backslash\{cone(e_{1}+me_{2},e_{2}),cone(e_{2})\}$.

We claim that, for each integer $r>1$, there is an unique variety
$\widetilde{Z}'_{r}$ with the two following properties: 1) the fan
$\widetilde{\Delta}'_{r}$ of $\widetilde{Z}'_{r}$ has support
$cone(e_{1},e_{1}+re_{2})$; 2) $\widetilde{Z}'_{r}$ is an open
subvariety of a toric variety  $\widetilde{Z}_{r}$ with nef
anticanonical bundle and birationally proper over $\mathbb{A}^{2}$.
In particular, the anticanonical divisor of $\widetilde{Z}'_{r}$ is
nef. The open subvariety $Z'_{r}$ of $Z_r$ whose fan is
$\Delta_{r}\backslash\{cone(e_{1}+re_{2},e_{2}),cone(e_{2})\}$,
satisfies these properties. So it is sufficient to prove the claim.

We show the claim for induction on $r$. We have already verified the
basis of induction. Let $\widetilde{Z}'_{r}$ be a variety that
satisfies the hypotheses of the claim and let $\sigma'$ be the
unique cone in $\widetilde{\Delta}_{r}'(2)$ which contains
$e_{1}+re_{2}$. Because of the inductive hypothesis  it is
sufficient to show that
$\sigma'=cone(e_{1}+re_{2},e_{1}+(r-1)e_{2})$.

Let $k$ be the function associated to the anticanonical bundle of a
fixed $\widetilde{Z}_{r}$. 
Let
$k_{\sigma'}=x_{1}e^{*}_{1}+x_{2}e^{*}_{2}$, then
$1=k_{\sigma'}(e_{1}+re_{2})=x_{1}+rx_{2}$ and
$x_{i}=k_{\sigma'}(e_{i})\leq1$ for each $i$, so the unique
possibilities for $k_{\sigma'} $ are $e^{*}_{1}$ and
$-(r-1)e^{*}_{1}+e^{*}_{2}$. Write
$\sigma'= cone(e_{1}+re_{2},v)$ with  $ v =c_{1}e_{1}+c_{2}e_{2}  $ primitive. If
$k_{\sigma'}=-(r-1)e^{*}_{1}+e^{*}_{2}$, then $c_{2}=(r-1)c_{1}+1$
because $(-(r-1)e^{*}_{1}+e^{*}_{2})(v)=1$. Because of the
non-singularity of $Z$ we have $c_{1}-1=\pm 1$, so there are two
possibilities: either $\sigma'=cone(e_{1}+re_{2},e_{2})$ or
$\sigma'=cone(e_{1}+re_{2},2e_{1}+(2r-1)e_{2})$. We exclude
the first one because $e_{2}$ does not belong to
$|\widetilde{\Delta}'_{r}|$. We exclude also the second one because
$k_{cone(e_{1},e_{1}+e_{2})}
(2e_{1}+(2r-1)e_{2})=
2>k(v) $. If $k_{\sigma'}=e^{*}_{1}$, then $c_{1}=k_{\sigma'}(v)=1$.
Because of the smoothness of $Z$ we have $c_{2}-r=\pm 1$. Again, we exclude
$e_{1}+(r+1)e_{2}$  because it does not belong to    $|\widetilde{\Delta}'_{r}|$.
Thus $\sigma'=cone(e_{1}+re_{2},e_{1}+(r-1)e_{2})$.
$\square$

Suppose $\grD_X[p]=\{e_1,e_1+e_2,e_1+2e_2,...,e_1+re_2,e_2\}$. Then $-K_X$ is nef if and only  if $-2\rho +2\rho_0 +e^*_{1}$ and $-2\rho +2\rho_0 -(r-1)e^*_{1}+e^*_{2}$ are antidominant. In such a case, if $\grt$ is decomposable, then the sum of these weights is regular, so $-K_X$ is big.  If $r=1$ and the previous weights are regular, then $-K_X$ is ample.
Remark that if $r=1$ and $\theta$ is decomposable, such weights are always regular.

\

III) If $X$ is not toroidal and $\grt$ is indecomposable, let $X'$ be as in Lemma \ref{riduzione a fano}. By Lemma \ref{lem fano indecomposable non toroidal} and Theorem \ref{fano rank2}, $R_{G,\theta}=A_2$, $H=G^\grt$ and  $\grD_{X'}[p]$ contains properly $\{\gra^\cech_1,-\gra^\cech_1-\gra^\cech_2,-3\gro^\cech_1\}$, up to reindexing.
By Lemma \ref{toroidal qfano rango2}  we have two possibilities for  $\grD_{X'}[p]$: i) $\{\gra^\cech_1,-\gro^\cech_1-\gro^\cech_2,-4\gro^\cech_1-\gro^\cech_2,
-5\gro^\cech_1- 2\gro^\cech_2,...,
-(r+3)\gro^\cech_1-r\gro^\cech_2,-3\gro^\cech_1\}$ and ii)  $\{\gra^\cech_1,-\gro^\cech_1-\gro^\cech_2,-4\gro^\cech_1- \gro^\cech_2,
-7\gro^\cech_1-\gro^\cech_2,...,
-(3r+1)\gro^\cech_1-\gro^\cech_2,-3\gro^\cech_1\}$.

Recall that $a_i=(2\rho-2\rho_0)(\gra_i^\cech)$. In the first case $k$ is linear over   $cone(\gra^\cech_1,-\gro^\cech_1-\gro^\cech_2)$,   $\grs:=cone(-\gro^\cech_1-\gro^\cech_2,-(r+3)\gro^\cech_1-r\gro^\cech_2)$ and  $cone(-(r+3)\gro^\cech_1-r\gro^\cech_2,-3\gro_1^\cech)$. The unique non-trivial condition is $k_{\grs}(   \gra^\cech_1)=r  \leq a_1$.
In the second case $k$ is linear over   $cone(\gro^\cech_1,-\gro^\cech_1-\gro^\cech_2)$,   $cone(-\gro^\cech_1-\gro^\cech_2,-(3r+1)\gro^\cech_1- \gro^\cech_2)$ and  $\grs:=cone(-(3r+1)\gro^\cech_1- \gro^\cech_2,-3\gro_1^\cech)$. The unique non-trivial condition is   $k_{\grs}( \gra^\cech_2)=2r-1\leq a_2$ (or equivalently $\widetilde{k}_{\grs}( \gra^\cech_2 ) \leq 0$). If $r=1$, then $\grt=-id$ over $\chi_*(T)$ so that  $-K_X$ is not ample.

IV) Finally, suppose $X$ non toroidal and $\grt$ decomposable. Suppose also $X_0$ smooth; the case $X_0$ singular is very similar.  By the proof of Theorem \ref{fano rank2}, $\rho(D(X))$ cannot be $\rho(D(G/H))$.
If   $\rho(D(X))= \alpha_{1}^{\vee}$, then, by the local factoriality of $X$, $e_{1}=-\alpha_{1}^{\vee}$ and there is $r\in\mZ^{>0}$ such that
$\sigma:=cone(\alpha_{1}^{\vee},v'=-r\alpha_{1}^{\vee}+e_{2})\in \Delta$. We apply Lemma    \ref{riduzione a fano} to the maximal open toroidal subvariety of $X$ (whose colored fan has support $cone(-\gra^\cech_1,v')$). Then, by Lemma  \ref{toroidal qfano rango2},
$ \Delta[p]$ has to be
$\{-\alpha_{1}^{\vee},v',-\alpha_{1}^{\vee}+v',-\alpha_{1}^{\vee}+2v',...,-\alpha_{1}^{\vee}+sv',\alpha_{1}^{\vee}\}$,
$\{-\alpha_{1}^{\vee},v',-\alpha_{1}^{\vee}+v',-2\alpha_{1}^{\vee}+v',...,-s\alpha_{1}^{\vee}+v',\alpha_{1}^{\vee}\}$ or $\{-\alpha_{1}^{\vee},v',\alpha_{1}^{\vee}\}$.

In the first case, 
$k$ is linear over $
cone(\gra^\cech_1,v')$, $\grs_1:=cone(v',-\alpha_{1}^{\vee}+sv')$ and $\grs_2:=cone(-\alpha_{1}^{\vee}+sv',-\gra^\cech_1)$.
The non-trivial conditions are $\tilde{k}_{\grs_1}(\gra^\cech_1) \leq0$ and $\tilde{k}_{\grs_2}( \gra^\cech_2)\leq0$ (or equivalently  $s\leq \bar{m}_1+1$ and $r\leq \bar{m}_2$). If $s=1$, then $r>m_2$ because $-K_X$ is not ample.
%
%
In the second case we can suppose $s\geq2$; 
$k$ is linear over  $
cone(\gra^\cech_1,v')$, $\grs_1:=cone(v',-s\alpha_{1}^{\vee}+v')$ and $\grs_2:=cone(-s\alpha_{1}^{\vee}+v',-\gra^\cech_1)$.
The unique non-trivial condition is $\tilde{k}_{\grs_2}( \gra^\cech_2)\leq 0$, or, equivalently, $r+s-1\leq \bar{m}_2$.
In the first two cases $\tilde{k}_{\grs_1}(\gra^\cech_2)<0$ and
$\tilde{k}_{\grs_2}(\gra^\cech_1)<0$, so $-\widetilde{K}_X $ is big. In the last case, we proceed as in the proof of Theorem \ref{fano rank2}, obtaing $m_2+1<r\leq \bar{m}_2+1$. Moreover, $\tilde{k}_{cone(\gra_1^\cech,v')}(\gra^\cech_2)<0$ and
$\tilde{k}_{cone(-\gra_1^\cech,v')}(\gra^\cech_1)<0$, so $-K_X$ is big.

We can study the smoothness of all the previous varieties as in the proof of Theorem \ref{fano rank2}. Remark that if
$H=G^\grt$, $\grt$ is decomposable,   $\sharp\rho^{-1}(\gra^\cech_1)=2$  and $(C(Y),D(Y))= (cone(\gra^\cech_1,-r\gra^\cech_1-\gra^\cech_2),F)$, then $Y$ is smooth if and only if $\sharp F=1$ and $\gra_1$ is not exceptional (but $Y$ cannot be an open subvarieties of a Fano variety).
The symmetric varieties in the statement are all   projective. Indeed, if $\Delta_X[p]= \{e_{1},e_{1}+e_{2},e_{1}+2e_{2},...,
e_{1}+(s-1)e_{2},e_{1}+se_{2},e_{2})\}$, then the following divisor is ample: 
$-pK_X+\sum_{i=1}^{s}i^2D_{e_1+ie_2}+ m\sum_{D(G/H)\setminus D(X)}D$ with $p,m >>0$. Indeed the piecewise linear function associated to $\sum_{i=1}^{s}i^2D_{e_1+ie_2}$ is strictly convex on
$\{C\in \grD_X: C\subset cone(e_1,e_1+s_2)\}$. Moreover, $k_X$ is strictly convex on the fan with maximal cones $cone(e_1,e_1+s_2)$ and $cone( e_1+s_2,e_2)$. The other cases are similar to this one.

\begin{figure}\begin{changemargin}{-2cm}{1cm}
\begin{tabular}{ |l|l|l|}\hline

 $G/H$ &$\theta$& $\Delta[p]$\\\hline \hline

$PSL_{3}$&$(A_{2})$& $\nexists$\\\hline

$SL_{3}$&$(A_{2})$&  $\{\gra^\cech_i, -\gro^\cech_i-\gro^\cech_{i^c},-4\gro^\cech_i- \gro^\cech_{i^c},-5\gro^\cech_i-2\gro^\cech_{i^c}, -3\gro^\cech_i \}$, $i=1,2$\\\hline

$SO_{5}$&($B_{2}$)& $\{-\gro^\cech_{1},-\gro^\cech_{2},-\gro^\cech_{1}-\gro^\cech_{2}\}$\\\hline

$Spin_{5}$&($B_{2}$)& $ \{-2\gro^\cech_{1},-\gro^\cech_{2},-2\gro^\cech_{1}-\gro^\cech_{2},
-2\gro^\cech_{1}-2\gro^\cech_{2}\} $  \\\hline

$ G_{2}$& $(G_2)$&$\nexists$ \\ \hline

$SL_{3}/N (SO_{3})$ &$AI$  & $\{-\gro_{1}^{\vee},-\gro_{2}^{\vee},-\gro_{1}^{\vee}-\gro_{2}^{\vee}\}$\\\hline

$ SL_{3}/SO_{3}$ &$AI$ & $\{ \alpha_{i}^{\vee},-3\gro_{i}^{\vee},-\gro_{i}^{\vee}-\gro_{i^c}^{\vee},-4\gro^\cech_i-\gro^\cech_{i^c}\}$, $i=1,2$\\ \hline

$SL_{6}/N (Sp_{6})$ &$AII$& $ \{- \gro^\cech_{1},-\gro^\cech_{2},- \gro^\cech_{i}-\gro^\cech_{i^c},
-2\gro^\cech_{i}- \gro^\cech_{i^c}\} $, $i=1,2$ \\\hline

$SL_{6}/Sp_{6}$ &$AII$&    $\{\gra^\cech_i, -\gro^\cech_i-\gro^\cech_{i^c},
...,-(r+3)\gro^\cech_i-r\gro^\cech_{i^c},-3\gro^\cech_i \}$, $i=1,2$; $r=2,3,4 $
 \\\cline{3-3}
 & &  $\{\gra^\cech_i, -\gro^\cech_i-\gro^\cech_{i^c},-4\gro^\cech_i- \gro^\cech_{i^c},-7\gro^\cech_i- \gro^\cech_{i^c},-3\gro^\cech_i \} $, $i=1,2$\\\cline{3-3}
\hline

$ SL_{n+1}/S(GL_{2}\times GL_{n-1})$ &$AIII$ & $\{-\gro_{1}^{\vee},-\gro_{2}^{\vee},-\gro_{1}^{\vee}-\gro_{2}^{\vee},
...,- r \gro_{1}^{\vee}-\gro_{2}^{\vee}\}$, $2\leq r\leq n-2$\\\cline{3-3} $n\geq4$&
&$\{\gra^\cech_2,-\gro^\cech_2\}$\\\hline

$ SL_{4}/N (S(GL_{2}\times GL_{2}))$ &$AIII$&$\{-\gro_{1}^{\vee},-\gro_{2}^{\vee},-\gro_{1}^{\vee}-\gro_{2}^{\vee}\}$\\\hline

$ SL_{4}/S(GL_{2}\times GL_{2})$ &$AIII$& $\{-\gro_{1}^{\vee},-2\gro_{2}^{\vee},-\gro_{1}^{\vee}-2\gro_{2}^{\vee}\}$\\\cline{3-3}&&$\{\gra^\cech_2,-2\gro^\cech_2\}$\\\hline

$ SO_{2n+1}/S(O_2\times O_{2n-1}) $&$BI$& $\{-\gro_{1}^{\vee},-\gro_{2}^{\vee},-\gro_{1}^{\vee}-\gro_{2}^{\vee},
...,-r\gro_{1}^{\vee}-\gro_{2}^{\vee}\}$,  $r\leq n-2$; $n\geq3$\\
\hline

$ SO_{5}/S(O_2\times O_{3}) $&$BI$&$\nexists$  \\\hline

$ SO_{2n+1}/  $&$BI$& $\{-2\gro_{1}^{\vee},-\gro_{2}^{\vee},-2\gro_{1}^{\vee}-\gro_{2}^{\vee},
...,-2r\gro_{1}^{\vee}-\gro_{2}^{\vee}\}$,  $r\leq n-1$
\\\cline{3-3}$(SO_2\times SO_{2n-1})$
&&$\{\gra^\cech_1,-\gro^\cech_1\}$\\\hline

$ Sp_{2n}/(Sp_{4}\times Sp_{2n-4}),$&$CII$& $\{-\gro_{1}^{\vee},-\gro_{2}^{\vee},-\gro_{1}^{\vee}-\gro_{2}^{\vee},
...,-r\gro_{1}^{\vee}-\gro_{2}^{\vee}\}$, $2\leq r\leq n-5$
\\\cline{3-3}$n>4$
&& $ \{- \gro^\cech_{1},-\gro^\cech_{2},- \gro^\cech_{1}-\gro^\cech_{2},
- \gro^\cech_{1}- 2\gro^\cech_{2}\} $\\\hline

$ Sp_{8}/N(Sp_{4}\times Sp_{4})$ &$CII$&$\{-\gro_{1}^{\vee},-\gro_{2}^{\vee},-\gro_{1}^{\vee}-\gro_{2}^{\vee},-2\gro_{1}^{\vee}-\gro_{2}^{\vee}\}$\\\cline{3-3}
&&$\{-\gro_{1}^{\vee},-\gro_{2}^{\vee},-\gro_{1}^{\vee}-\gro_{2}^{\vee},-\gro_{1}^{\vee}-2\gro_{2}^{\vee}\}$\\\hline

$Sp_{8}/(Sp_{4}\times Sp_{4})$ & $CII$& $\{-\gro_{1}^{\vee},-2\gro_{2}^{\vee},-\gro_{1}^{\vee}-2\gro_{2}^{\vee},
...,-r\gro_{1}^{\vee}-2 \gro_{2}^{\vee}\}$, $2\leq r\leq3$ \\\cline{3-3}&&$ \{- \gro^\cech_{1},-2\gro^\cech_{2},- \gro^\cech_{1}-2\gro^\cech_{2},
- \gro^\cech_{1}- 4\gro^\cech_{2}\} $\\\hline

$SO_{2n}/S(O_{2}\times O_{2n-2})$&$DI$&  $\{-\gro_{1}^{\vee},-\gro_{2}^{\vee},-\gro_{1}^{\vee}-\gro_{2}^{\vee},...,-r\gro_{1}^{\vee}-\gro_{2}^{\vee}\}$, $r\leq n-2$\\\hline

$SO_{2n}/(SO_{2}\times SO_{2n-2})$ &$DI$&  $\{-2\gro_{1}^{\vee},-\gro_{2}^{\vee},-2\gro_{1}^{\vee}-\gro_{2}^{\vee},...,-2r\gro_{1}^{\vee}-\gro_{2}^{\vee}\}$,  $r\leq n-2$ \\\cline{3-3}&&$\{\gra^\cech_1,-2\gro^\cech_1\}$\\\hline

$ SO_{8}/N (SO_{4})$ &$DIII$ &  $\{-\gro_{1}^{\vee},-\gro_{2}^{\vee},-\gro_{1}^{\vee}-\gro_{2}^{\vee}\}$\\\cline{3-3}
&&$\{-\gro_{1}^{\vee},-\gro_{2}^{\vee},-\gro_{1}^{\vee}-\gro_{2}^{\vee},-\gro_{1}^{\vee}-2\gro_{2}^{\vee}\}$\\\hline

$ SO_{8}/SO_{4}$ &  $DIII$ &  $\{-\gro_{1}^{\vee},-2\gro_{2}^{\vee},...,-\gro_{1}^{\vee}-2r\gro_{2}^{\vee}\}$, $r=1,2$\\\cline{3-3}
&&$\{\gra^\cech_2,-2\gro^\cech_2\}$\\\hline

$SO_{10}/SO_{5}$ &$DIII$ & $\{-\gro_{1}^{\vee},-\gro_{2}^{\vee},-\gro_{1}^{\vee}-\gro_{2}^{\vee},...
-r\gro_{1}^{\vee}-\gro_{2}^{\vee}\}$, $r=2,3$ \\\cline{3-3}
 &  & $\{-\gro_{1}^{\vee},-\gro_{2}^{\vee},-\gro_{1}^{\vee}-\gro_{2}^{\vee},
- \gro_{1}^{\vee}-2\gro_{2}^{\vee}\}$
\\\cline{3-3}&&$\{\gra^\cech_2,-\gro^\cech_2\}$ \\\hline

$E_{6}/(D_5\times \mC^*)$&$EIII$& $\{-\gro_{1}^{\vee},-\gro_{2}^{\vee},-\gro_{1}^{\vee}-\gro_{2}^{\vee},...,-r\gro_1^{\vee}-\gro_2^{\vee}\}$,  $2\leq r\leq5$\\\cline{3-3}
 & & $\{-\gro_{1}^{\vee},-\gro_{2}^{\vee},-\gro_{1}^{\vee}-\gro_{2}^{\vee}, ...,- \gro_{1}^{\vee}-r\gro_{2}^{\vee}\}$, $2\leq r\leq3$
\\\cline{3-3}&&$\{\gra^\cech_2,-\gro^\cech_2\}$\\\hline

$E_{6}/N (F_{4})$&$EIV$&$\{-\gro_{1}^{\vee},-\gro_{2}^{\vee},-\gro_{i}^{\vee}-\gro_{i^c}^{\vee},...-r\gro_{i}^{\vee}-\gro_{i^c}^{\vee}\}$, $i=1,2$; $r=2,3,4$\\\hline

$ E_{6}/F_{4}$&$EIV$ &    $\{\gra^\cech_i, -\gro^\cech_i-\gro^\cech_{i^c},...,-(r+3)\gro^\cech_i-r\gro^\cech_{i^c},-3\gro^\cech_i \}$, $i=1,2$; $r=2,...,8$
\\\cline{3-3}
 & &  $\{\gra^\cech_i, -\gro^\cech_i-\gro^\cech_{i^c},..., -(3r+1)\gro^\cech_i- \gro^\cech_{i^c},-3\gro^\cech_i \}$, $i=1,2$; $r=2,3,4$\\\hline

$ G_{2}/(A_{1}\times A_{1})$&$G$&  $\{-\gro_{1}^{\vee},-\gro_{2}^{\vee} \}$\\
\hline

\end{tabular}\end{changemargin}
\caption{quasi-Fano indecomposable symmetric varieties with rank 2}
\label{fig qFano rank2 indecomposable}
\end{figure}


%
%
%
%
%
%
%
%
%
%
%
%
%
%
%

\subsection{Symmetric Fano varieties with $G$ reductive}\label{sez rang2+G red}

In this section we consider the 2-rank locally factorial Fano symmetric varieties over which acts a group $G$ which is only reductive. If $G$ is
a torus, then $X$ is the projective space. So, we can suppose
$G=G'\times \mathbb{C}^{*}$ with $H\cap
\mathbb{C}^{*}=(\mathbb{C}^{*})^{\theta}=\{\pm id\}$. Write
$\chi_{*}(\mathbb{C}^{*}/\{\pm id\})=\mathbb{Z}f$.

If $R_{G,\theta}=BC_{1}$, then   $H=G^\grt$ and
$\chi_{*}(S)=\mathbb{Z}f\oplus\mathbb{Z}\alpha^{\vee} $. Instead,
if $R_{G,\theta}=A_{1}$ there are three possibilities: 1)
$\chi_{*}(S)=\mathbb{Z}f\oplus\mathbb{Z}\alpha^{\vee}$ and
$H=G^{\theta}$; 2)
$\chi_{*}(S)=\mathbb{Z}f\oplus\mathbb{Z}\frac{1}{2}\alpha^{\vee}$
and $H=N_{G'}(G^{\theta})\times \{\pm1\}$;  3)
$\chi_{*}(S)=\mathbb{Z}f\oplus\mathbb{Z}\frac{\alpha^{\vee}+f}{2}$.
In the last case $H$ is
generated by $G^{\theta}$ and by $n_1n_2$, where $n_1\in
N_{G'}((G')^{\theta})\setminus (G')^{\theta}$ and  $n_2\in\mathbb{C}^{*}$ has  order four; in particular $[G^{\theta}:H]=2$. Let $e$ be the primitive positive multiple of $-\gra^\cech$ and let $\{e^*,f^{*}\}$ be
the dual basis of $\{e ,f\}$.

\begin{thm}\label{thm X fano G reductive rango2} Let $G/H$ be a   symmetric space of rank two, such  that $G$ is neither   semisimple nor abelian. As  before, write
$\psi (r)= -2\rho +2\rho_0 -re^* $ and $m_1:=max\{r: \psi(r)<0\}$. The   Fano locally factorial compactifications of $G/H$ are the following ones:
\begin{enumerate}

\item $\Delta[p]=\{f,-f,e+rf\}$ if  $\chi_{*}(S)=\mathbb{Z}e\oplus\mathbb{Z}f $, $r\in\mathbb{Z}$,  $r\leq m_1+1$ and $-r\leq m_1+1$;
\item $\Delta[p]=\{f,-f,e+rf, e+(r+1)f\}$    if $\chi_{*}(S)=\mathbb{Z}e\oplus\mathbb{Z}f $, $r\in\mathbb{Z}$, $r\leq m_1$ and $-r\leq m_1+1$;
\item $\Delta[p]= \{\alpha^{\vee},-f,-\alpha^{\vee}+f\}$ (and
$\Delta[p]=\{\alpha^{\vee},f,-\alpha^{\vee}-f\}$) if  $\chi_{*}(S)=\mathbb{Z}\alpha^{\vee}\oplus\mathbb{Z}f $ and $\sharp D(G/H)=1$;
\item $\Delta[p]=\{\alpha^{\vee},-f,-\alpha^{\vee},-\alpha^{\vee}+f\}$ (and $\Delta[p]=\{\alpha^{\vee},f,-\alpha^{\vee},-\alpha^{\vee}-f\}$) if $\chi_{*}(S)=\mathbb{Z}\alpha^{\vee}\oplus\mathbb{Z}f $ and $\sharp D(G/H)=1$;
\item $\Delta[p]=\{f,-f, -\frac{1}{2}\alpha^{\vee}+\frac{2r+1}{2}f\}$   if  $\chi_{*}(S)=\mathbb{Z}f\oplus\mathbb{Z} (\frac{\alpha^{\vee}+f}{2} ) $,  $r\in\mathbb{Z}$, $r\leq \frac{m_1+1}{2}$ and $-r\leq \frac{m_1+3}{2}$;
\item $\Delta[p]=\{f,-f, -\frac{1}{2}\alpha^{\vee}+\frac{2r+1}{2}f,-\frac{1}{2}\alpha^{\vee}+\frac{2r+3}{2}f\}$    if $\chi_{*}(S)=\mathbb{Z}f\oplus\mathbb{Z}(\frac{\alpha^{\vee}+f}{2}) $, $r\in\mathbb{Z}$, $r\leq \frac{m_1-1}{2}$ and $-r\leq \frac{m_1+3}{2}$;
\item $\Delta[p]=\{
\alpha^{\vee}, -\frac{1}{2}\alpha^{\vee}+
\frac{1}{2}f,-\frac{1}{2}\alpha^{\vee}- \frac{1}{2}f,-\alpha^{\vee}\}$ if $\chi_{*}(S)=\mathbb{Z}f\oplus\mathbb{Z}(\frac{\alpha^{\vee}+f}{2})  $;
\item $\Delta[p]=\{
\alpha^{\vee}, -\frac{1}{2}\alpha^{\vee} + \frac{1}{2}f,-\frac{
1}{2}\alpha^{\vee}- \frac{1}{2}f\}$ if  $\chi_{*}(S)=\mathbb{Z}f\oplus\mathbb{Z}(\frac{\alpha^{\vee}+f}{2}) $;

\item $\Delta[p]=\{\alpha^{\vee},-f, -\frac{1}{2}\alpha^{\vee}+ \frac{1}{2}f\}$  (and $\Delta[p]=\{\alpha^{\vee},f, -\frac{ 1}{2}\alpha^{\vee}- \frac{1}{2}f\}$) if $\chi_{*}(S)=\mathbb{Z}f\oplus\mathbb{Z}(\frac{\alpha^{\vee}+f}{2}) $;
\item $\Delta[p]=\{\alpha^{\vee},-f, -\frac{1}{2}\alpha^{\vee}+ \frac{1}{2}f,-\frac{1}{2}\alpha^{\vee}- \frac{1}{2}f\}$  (and $\Delta[p]=\{\alpha^{\vee},f,-\frac{1}{2}\alpha^{\vee}+ \frac{1}{2}f,-\frac{ 1}{2}\alpha^{\vee}- \frac{1}{2}f\}$) if $\chi_{*}(S)=\mathbb{Z}f\oplus\mathbb{Z}(\frac{\alpha^{\vee}+f}{2}) $.
\end{enumerate}
The only singular varieties are the ones in the cases 3) and 4).
\end{thm}

Observe that $X$ cannot be simple because the valuation cone is not strictly convex. We begin with a lemma similar to   Lemma \ref{numero 1 coni X fano rango 2}.

\begin{lem}\label{numero 1 coni G reductive} Let  $X$ be a Fano locally factorial symmetric variety with   $G$ as before. If $X$ is toroidal,  there are at most four colored 1-cones. Otherwise, there are at most three colored 1-cones.
\end{lem}

{\em Proof of Lemma \ref{numero 1 coni G reductive}.} First suppose $X$ toroidal. Then $k(\pm f)=1$. Let
$\sigma\in \Delta(2)$ be a cone which does not contain neither  $f$ nor $-f$, then
$k_{\sigma}(\pm f)\leq0$. So $k_{\sigma}$ is a multiple of
$e^{*}$. If  there is another cone
$\sigma'\in \Delta(2)$ with the same properties (and such that $dim\,(\sigma\cap \sigma')=1$), then   $k_{\sigma}=k_{\sigma'}$.
If $X$ is not toroidal, we can study the maximal open toroidal subvariety $X'$ of $X$ as in Lemma \ref{numero 1 coni X fano rango 2} (because $|\calF(X')|$ is strictly convex). $\square$


{\em Proof of Theorem  \ref{thm X fano G reductive rango2}.} We have to request, as in Theorem \ref{fano rank2}, that $\sharp\rho(D(X))$ $=\sharp\rho^{-1}(\rho(D(X)))$, but in this case $\sharp\rho(D(X))\leq1$.  If $X$ is not toroidal, then $e=-\gra^\cech$, because there is an appropriate $cone(\gra^\cech,v)$ in $\grD(2)$.
Suppose first that $\chi_{*}(S)=$ $\mathbb{Z}e$ $\oplus $ $\mathbb{Z}f$.  We have to consider the following cases: 1) If $\Delta[p]=\{f,- f,v\}$, then $v=e+rf$ because $\{ f,v\}$ is a basis of $\chi_{*}(S)$.  We have to impose that $\widetilde{k}_{cone(f,e+rf)}(\gra^\cech)=\psi(r-1)(\gra^\cech)<0$ and that $\widetilde{k}_{cone(-f,e+rf)}(\gra^\cech)=\psi(-r-1)(\gra^\cech)<0$.
2) Suppose $\Delta[p]=\{f,-f,v_{1},v_{2}\}$. As before
$v_{i}=e+x_{i}f$; moreover  $x_{1}=x_{2}\pm1$
because $\{v_{1},v_{2}\}$ is a basis of $\chi_{*}(S)$. Suppose $x_1=x_2+1$. We have to impose  $\widetilde{k}_{cone(f,e+(x_2+1)f)}(\gra^\cech)=\psi(x_2 )(\gra^\cech)<0$ and $\widetilde{k}_{cone(-f,e+x_2f)}(\gra^\cech)=\psi(-x_2-1)(\gra^\cech)<0$.


3) Suppose $\Delta[p]=\{\alpha^{\vee},-f,v\}$. Then
$v=-\gra^\cech  +rf=
-m\alpha^{\vee}\pm f$ because $X$ is smooth. But $-f$ is not contained in
$cone(\alpha^{\vee},v)$, thus
$v=-\alpha^{\vee}+f$.
4) Suppose $\Delta[p]=\{\alpha^{\vee},-f,v_{1},v_{2}\}$ and
$\Delta(2)=\{cone(\alpha^{\vee},v_{1}), \sigma: =cone(v_{1},$ $v_{2}),$ $
cone(-f,v_{2})\}$. Then    $v_{1}=-m\alpha^{\vee}
+ f$ and $v_{2}=-\gra^\cech+rf$.   Observe that
$\mathbb{R}^{\geq0}v_{1}=\mathbb{R}^{\geq0}(-\alpha^{\vee}+ \frac{1}{ m}f)$, so
$r\leq0$. Furthermore $v_{2}-rv_{1}=-(1-r m)\gra^\cech=\pm \gra^\cech$, so
$ mr$ is 0 or 2. 
Thus   $v_{2} =-\alpha^{\vee}$. Moreover
$h_{\sigma}(-f)=
m-1<1$, so  $v_{1}=-\alpha^{\vee}+f$.


5) Suppose $\Delta[p]=\{\alpha^{\vee},v_{1},v_{2}\}$ and
$\Delta(2)=\{cone(\alpha^{\vee},v_{1}),$ $ cone(v_{1},$ $v_{2}),$ $
cone(\alpha^{\vee},v_{2})\}$. Thus $v_{1}=-r\alpha^{\vee}+f$ and
$v_{2}=-m\alpha^{\vee}-f$ with $r,m>0$. Moreover,
$v_{1}+v_{2}=-(r+m)\alpha^{\vee}=\pm\alpha^{\vee}$, a contradiction.
6) Suppose $\Delta[p]=\{\alpha^{\vee},v_{1},v_{2}, v_3\}$ and
$\Delta(2)=\{cone(\alpha^{\vee},v_{1}), \grs:=cone(v_{1},v_{2}), cone(v_{2},v_{3}),$ $
cone$ $(\alpha^{\vee},$ $v_{3})\}$. Write $v_2=-x\gra^\cech+yf$. By the smoothness of $X$, $v_{1}=-r\alpha^{\vee}+f$,
$v_{3}=-m\alpha^{\vee}-f$, $x=-my\pm1$ and $x=+ry\pm1$. The last two conditions, plus $x>0$ imply that 
$v_2=-\gra^\cech$.  Moreover $k_{\grs}(v_3)=
m+r-1<1$,  a contradiction.

Now,  suppose
$\chi_{*}(S)=\mathbb{Z}\frac{f+\alpha^{\vee}}{2}\oplus\mathbb{Z}f$.
1) The toroidal case can be studied as before.
2) Suppose that $\Delta(2)$ contains two cones
$\sigma_{\pm}:=cone(v_{\pm},\alpha^{\vee})$.  Let $u=\frac{1}{2}\gra^\cech+\frac{1}{2}f$.
We have
$v_{+}=-\frac{2m +1}{2}\alpha^{\vee}+ \frac{1}{2}f=-(2m+1)u+(m+1)f$ and
$v_{-}=-\frac{2r +1}{2}\alpha^{\vee}- \frac{1}{2}f=-(2r+1)u+rf$ with $m,r\geq0$.
First, suppose that there is another
$v=-xu+yf$ in $\Delta[p]$; here $x>0$. By the smoothness, we have $x(m+1)=y(2m+1)\pm1$ and $xr=y(2r+1)\pm1$. Thus
$x[(m+1)(2r+1)-r(2m+1)]= \pm1\mp1$, so $x(m+r+1)=2$, $x(m+1)=y(2m+1)+1$ and $xr=y(2r+1)-1$. If $x=1$, then the previous three equations are not compatible. Instead, if $x=2$ then 
$\grD[p]=\{\gra^\cech,-\frac{1}{2}\gra^\cech+\frac{1}{2}f,- \gra^\cech,-\frac{1}{2}\gra^\cech-\frac{1}{2}f\}$.
Next, suppose that $\Delta[p]=\{\alpha^{\vee},v_{+}, v_-\}$; in particular
$\{v_{+},v_{-}\}$ is a basis. Thus
$v_{+}+v_{-}=-(m+r +1)\alpha^{\vee}=\pm\alpha^{\vee}$, so
  $\Delta[p]=\{\alpha^{\vee},-\frac{1}{2}\alpha^{\vee}+ \frac{1}{2}f,
-\frac{1}{2}\alpha^{\vee}-\frac{1}{2}f \}$.


3) Finally, suppose that $\Delta[p]$ contains $\alpha^{\vee}$ and
$-f$. Then $\Delta$ contains $\sigma=cone(\alpha^{\vee},v_{1})$ and
$ cone(-f,v_{2})$, with
$v_{1}=-\frac{2m+1}{2}\alpha^{\vee}+\frac{1}{2}f$, $m\geq0$ and
$v_{2}=-\frac{1}{2}\alpha^{\vee}-\frac{2r+1}{2}f$.
If $v_1=v_2$, then $\grD[p]=\{\gra^\cech,-f,-\frac{ 1}{2}\alpha^{\vee}+\frac{1}{2}f \}$. Otherwise, by Lemma \ref{numero 1 coni G reductive}, $\grD$ contains $cone(v_1,v_2)$.
Furthermore, $(2r+1)v_{1}+v_{2}=
\pm\alpha^{\vee}$,
hence $(2r+1)(2m+1)$ is 1 or -3. But $r\geq0$ because $v_2\notin\grs$, so  $v_{1}=-\frac{1}{2}\alpha^{\vee}+f$
and $v_{2}=-\frac{1}{2}\alpha^{\vee}-f$.

\end{document}